\documentclass[a4paper,num-refs]{article}

\pdfoutput=1

\usepackage{IEK10} 
\usepackage{natbib}

\makeatother
\usepackage{amssymb}
\usepackage{xcolor}
\usepackage{mathtools, cuted}
\usepackage{tabularx}
\usepackage{eurosym}
\usepackage[utf8]{inputenc}
\usepackage{multicol}
\usepackage{multirow}

\definecolor{JuelichBlue}{RGB}{2, 61, 107}
\definecolor{JB}{RGB}{2, 61, 107}
\definecolor{JuelichLightBlue}{RGB}{173, 189, 227}
\definecolor{JLB}{RGB}{173, 189, 227}
\definecolor{JuelichRed}{RGB}{235, 95, 115}
\definecolor{JR}{RGB}{235, 95, 115}
\definecolor{G1}{RGB}{150, 150, 150}
\definecolor{JV}{RGB}{175,130,185}

\newcommand{\vect}[1]{\mathbf{#1}}

\renewcommand{\min}[1][]{
	\ifthenelse{\isempty{#1}}{\operatorname{min}}{\ensuremath{\underset{#1}{\text{min}\,}}}
}

\newcommand{\RNum}[1]{\uppercase\expandafter{\romannumeral #1\relax}}

\def\IEK10{
  Forschungszentrum Jülich GmbH,
  Institute of Energy and Climate Research,
  Energy Systems Engineering (IEK-10),
  Jülich 52425,
  Germany
}
\def\RWTH{
  RWTH Aachen University
  Aachen 52062,
  Germany
}
\def\ETH{
  ETH Zürich,
  Energy \& Process Systems Engineering,
  Zürich 8092,
  Switzerland
}

\newcommand{\mytitle}{ Demand Response for Flat Nonlinear MIMO Processes  using  \newline Dynamic Ramping Constraints}

\newcommand{\affil}{
  \begin{itemize}[leftmargin=3mm, itemsep=0mm]
    \item[$^a$]\IEK10
    \item[$^b$]\ETH
    \item[$^c$]\RWTH
  \end{itemize}
}

\def\firstAuthor{Florian Joseph Baader}
\newcommand{\myauthor}{\firstAuthor$^{a,b,c,*}$, Philipp Althaus$^{a,c}$, André Bardow$^{a,b}$, Manuel Dahmen$^{a}$}
\newcommand{\correspondingauthor}{F. J. Baader, \ETH \newline E-mail: fbaader@ethz.ch}
\author{\myauthor}

\usepackage[
  colorlinks,
  linkcolor=blue,
  citecolor=blue,
  urlcolor=blue,
  pdftitle={\mytitle},
  pdfauthor={\firstAuthor}
]{hyperref}
\usepackage[capitalise, nameinlink]{cleveref}
\crefname{table}{Tab.}{Tab.}
\crefname{equation}{equation}{equations}

\newcommand{\setpgfexternalcounter}[1]{
  \makeatletter%
  \pgfkeysgetvalue{/tikz/external/figure name}\myexternalname
  \expandafter\gdef\csname c@tikzext@no@\myexternalname\endcsname{#1}%
  \makeatother
}

\begin{document}

  \thispagestyle{firststyle}

  \begin{center}
    \begin{large}
    {\fontsize{12}{14} \selectfont
      \centering \textbf{\mytitle}}
    \end{large} \\
    \myauthor
  \end{center}

  \begin{footnotesize}
    \affil
  \end{footnotesize}

\begin{abstract}
Volatile electricity prices make demand response (DR) attractive for processes that can modulate their production rate. 
However, if nonlinear dynamic processes must be scheduled simultaneously with their local multi-energy system, the resulting scheduling optimization problems often cannot be solved in real time.
For single-input single-output processes, the problem can be simplified without sacrificing feasibility by dynamic ramping constraints that define a derivative of the production rate as the ramping degree of freedom.
In this work, we extend dynamic ramping constraints to flat multi-input multi-output processes by a coordinate transformation that gives the true nonlinear ramping limits.
Approximating these ramping limits by piecewise affine functions gives a mixed-integer linear formulation that guarantees feasible operation.
As a case study, dynamic ramping constraints are derived for a heated reactor-separator process that is subsequently scheduled simultaneously with its multi-energy system.
The dynamic ramping formulation bridges the gap between rigorous process models and simplified process representations for real-time scheduling.
\end{abstract}


\newpage

\section{Introduction}
Reducing greenhouse gas emissions requires increased renewable electricity production that, however, gives a fluctuating supply. 
This fluctuating supply can be compensated by consumers that react to time-varying electricity prices by shifting their demand in time in a so-called demand response (DR) \citep{Zhang.2016}.
DR can be attractive for energy-intensive production processes with the flexibility to modulate their production rate.
To exploit the DR potential, a scheduling optimization is needed, which typically determines operational set-points for a time horizon in the order of one day \citep{Baldea.2014}.
However, such a scheduling optimization is computationally challenging for nonlinear processes that exhibit scheduling-relevant dynamics. 
The scheduling optimization becomes even more difficult 
if processes do not only consume electricity but also heating or cooling as these processes need to be scheduled simultaneously with the local multi-energy supply system (often also referred to as utility system)  \citep{Leenders.2019b}.
Local multi-energy supply systems bring integer on/off decisions into the scheduling optimization, especially as they often consist of multiple redundant units \citep{Voll.2013}.
Thus, the simultaneous scheduling optimization problem usually is a nonlinear mixed-integer dynamic optimization (MIDO) problem \citep{baader2022simultaneous}.

Traditionally, such challenging scheduling MIDO problems are simplified by introducing static ramping constraints that define the rate of change of the production rate $\rho$ as degree of freedom $\nu$ and limit $\nu$ with constant ramping limits $\nu^{\text{min}}$, $\nu^{\text{max}}$, see e.g., \cite{Carrion.2006,SumitMitra.2012,RichardAdamson.2017,DanyanZhou.2017,Hoffmann.2021}:
    \begin{align}
        \label{eq:DR_MIMO_DRMIMO_ramp_constr}
        \dot{\rho} = \nu ~\text{with}~ \nu^{\text{min}}\leq\nu\leq\nu^{\text{max}}
    \end{align}
If additionally, the nonlinear energy demand of the production process is approximated as a piecewise-affine (PWA) function, the complete problem can be reformulated as a mixed-integer linear program (MILP) \citep{Schafer.2020}.

\begin{figure*}[h]
\centering
 \includegraphics[trim=0 0 205 4,clip]{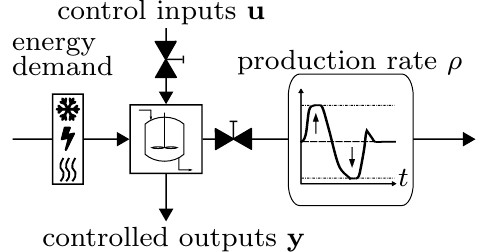}
  \caption{A MIMO process that can vary its production rate $\rho$ and thus its energy demand while additional control inputs $\vect{u}$ are available to control further process outputs variables $\vect{y}$.}
   \label{fig:DRMIMO_process_MIMO}
\end{figure*}

Traditional ramping constraints, however,  have two shortcomings: They are restricted to (i) first-order dynamics and (ii) constant ramping limits.
To tackle both shortcomings, we proposed high-order dynamic ramping constraints in our previous publication \citep{Baader.2021}.
These high-order dynamic ramping constraints limit the $\delta$-th derivative of the production rate $\rho$ and use limits $\nu^{\text{min}}$, $\nu^{\text{max}}$ which are functions of the production rate and its time derivatives:
\begin{align}
    \rho^{(\delta)} = \nu~\text{with}~ \nu^{\text{min}}\left(\rho,\dot{\rho},...,\rho^{(\delta-1)}\right) \leq \nu \leq \nu^{\text{max}}\left(\rho,\dot{\rho},...,\rho^{(\delta-1)} \right) 
\end{align}
Dynamic ramping constraints allow to represent energy-intensive processes better than static ones due to their ability to account for high-order dynamics and non-constant ramping limits that depend on the process state.
In \cite{Baader.2021}, we demonstrated that dynamic ramping constraints can achieve solutions of the scheduling optimization problem that are close to the solutions of the original nonlinear MIDO problem while  allowing to formulate MILPs that can be solved sufficiently fast for real-time scheduling.
Moreover, we presented a method to derive dynamic ramping constraints rigorously for the special case of single-input single-output (SISO) processes that are exact input-state linearizable. 
However, this case is quite restrictive, and the question remained how to derive dynamic ramping constraints for more general processes.
In particular, it was open how to apply the dynamic ramping method to multi-input multi-output (MIMO) processes that have a variable production rate $\rho$ while, at the same time, other output variables $\vect{y}$ such as temperatures and concentrations need to be controlled using a set of control inputs $\vect{u}$ (Figure \ref{fig:DRMIMO_process_MIMO}).

The present publication extends the dynamic ramping approach to MIMO processes and presents a method to derive dynamic ramping constraints rigorously for differentially flat MIMO processes.
In simple words, a MIMO process with $m$ inputs is differentially flat if an output vector $\boldsymbol{\xi}$ of the same size $m$ exists such that all  process states and inputs can be given as a function of the output $\boldsymbol{\xi}$ and its $\beta$ time derivatives $\dot{\boldsymbol{\xi}}$, ..., $\boldsymbol{\xi}^{(\beta)}$ \citep{Fliess.1995,Rothfu.1997,R.Rothfu.1996}. 
We make use of the fact that a flat nonlinear model can be transformed to a linear model \citep{Fliess.1995}.
However, constraints that are linear for the original model, e.g., bounds on inputs, are nonlinear for the transformed model.
In other words, a nonlinear model with linear constraints is transformed into a linear model with nonlinear constraints.
Approximating the nonlinear bounds to the safe side with piecewise-affine functions allows us to formulate a MILP, whose solution is guaranteed to be feasible also for the non-approximated version.
Note that for SISO processes, flatness is equivalent to exact input-state linearizability \citep{Adamy.2014}, which was the main assumption in our previous work \citep{Baader.2021}.
For MIMO processes, exact input-state linearizability with static state transformation is a special case of flatness \citep{Adamy.2014}.

The remaining paper is structured as follows: 
In Section \ref{sec:DRMIMO_Dyn_ramping_constraints}, the considered demand response optimization problem and the dynamic ramping constraints are introduced.
In Section \ref{sec:DRMIMO_method}, a method is presented to derive dynamic ramping constraints for flat MIMO processes.
In Section \ref{sec:DRMIMO_case_study}, a case study is presented.
Section \ref{sec:DRMIMO_discussion} provides further discussion and Section \ref{sec:DRMIMO_conlusion} concludes the work.

\section{Dynamic mixed-integer linear scheduling with ramping constraints}
\label{sec:DRMIMO_Dyn_ramping_constraints}
This section briefly introduces the simultaneous scheduling optimization problem (P) of flexible processes represented by dynamic ramping constraints and multi-energy supply systems.
This problem is a mixed-integer dynamic optimization (MIDO) problem with linear and piecewise affine (PWA) functions.
We discretize time through collocation on finite elements \citep{Biegler.2010} to convert the MIDO problem (P) to a MILP.
The size of the final MILP problem is proportional to the number of discretization points.
Note that all decision variables $\boldsymbol{\chi} = \left(\nu, (\vect{Q}_{\text{dem}}^{\text{process}})^T, \rho, S, \Phi_{\text{energy}},   (\vect{Q}^{\text{in}})^T, (\vect{Q}^{\text{out}})^T, (\Delta\vect{P})^T, \vect{z}_{\text{on}}^T \right)^T $, which are introduced in the following, are functions of time $t$. Still, we we do not state time dependency explicitly to ease readability.
The problem (P) without time discretization reads:
\begin{align}
    \label{eq:DR_MIMO_P1_objective} \tag{Pa}
    &\hspace{1.5cm}  \underset{\boldsymbol{\chi} \in \left[\boldsymbol{\chi}^l, \boldsymbol{\chi}^u\right]}{\text{min}}~~  \Phi_{\text{energy}}(t_f) \\
    \label{eq:DR_MIMO_P1_dynamic_ramping} \tag{Pb}
    \text{s.t. }~~& \text{Dynamic ramping constraints}\\
    \label{eq:DR_MIMO_P1_process_energy_dem} \tag{Pc}
    & \text{Process energy demand model}
    \\
    \label{eq:DR_MIMO_P1_storage} \tag{Pd}
    & \text{Product storage: } \dot{S} = \rho - \rho^{\text{nom}} ~~~\forall t \in \mathbb{T}
    \\ 
    \label{eq:DR_MIMO_P1_energy_costs} \tag{Pe}
    & \text{Energy costs: } \dot{\Phi}_{\text{energy}} = \sum_{e \in \mathbb{E}} p_e\left( \sum_{i \in \mathbb{C}_e^{\text{cons}}}  Q_{i}^{\text{in}} + \Delta P_e \right) ~~~\forall t \in \mathbb{T}
     \\
    \label{eq:DR_MIMO_P1_eff} \tag{Pf}
    & \text{Energy conversion: } Q_{i}^{\text{out}} = h_{Q_{i}^{\text{out}}}^{\text{PWA}}(Q_{i}^{\text{in}})  ~~~\forall i \in \mathbb{C},~\forall t \in \mathbb{T} 
    \\
    \label{eq:DR_MIMO_P1_part_load} \tag{Pg} 
    & \text{Minimum part-load : } z_i^{\text{on}} Q_i^{\text{min}} \leq Q_i^{\text{out}} \leq  z_i^{\text{on}} Q_i^{\text{max}} ,~~~\forall i \in \mathbb{C},~\forall t \in \mathbb{T} 
    \\ 
    \label{eq:DR_MIMO_P1_energy balance} \tag{Ph}
    & \text{Balances: } Q_{\text{dem},e}^{\text{process}} + Q_{\text{dem},e}^{\text{inflexible}} = \sum_{i \in \mathbb{C}_e^{\text{sup}}} Q_{i}^{\text{out}} +\Delta P_e, ~~~\forall e \in \mathbb{E},  ~\forall t \in \mathbb{T}
\end{align}
 The objective is to minimize the cumulative energy costs $\Phi_{\text{energy}}$ at final time $t_f$.
 The lower and upper bounds of the decision variables are $\boldsymbol{\chi}^l$ and $\boldsymbol{\chi}^u$.
In the following paragraph, we discuss the dynamic ramping constraints (\ref{eq:DR_MIMO_P1_dynamic_ramping}) and the process energy demand model (\ref{eq:DR_MIMO_P1_process_energy_dem}) in  detail.
The remaining constraints (\ref{eq:DR_MIMO_P1_storage}) to (\ref{eq:DR_MIMO_P1_energy balance}) are standard constraints \citep{Schafer.2020,SusanneSass.2020} and discussed in more detail in our previous publication \citep{Baader.2021}. Thus, we only briefly introduce the symbols here.
All decision variables are functions of time $t\in \mathbb{T}$ although not stated explicitly to ease readability. 
The product storage with level $S$ is filled by the production rate $\rho$ and emptied with the nominal production $\rho^{\text{nom}}$ (\ref{eq:DR_MIMO_P1_storage}). 
The rate of change of the energy costs $\Phi_{\text{energy}}$ is the price of an energy form $e$, $p_e$ times the input power $Q_{i}^{\text{in}}$ of energy conversion units consuming $e$, $\mathbb{C}_e$, and the power exchanged with the grid $\Delta P_e$. 
The set $\mathbb{E}$ in (\ref{eq:DR_MIMO_P1_energy_costs}) covers all considered energy forms.
For the energy conversion of a components $i$ in the set of components $\mathbb{C}$, the output power $Q_{i}^{\text{out}}$ is given as piecewise affine function $h_{Q_{i}^{\text{out}}}^{\text{PWA}}$ of the input power $Q_{i}^{\text{in}}$ (\ref{eq:DR_MIMO_P1_eff}).
Additionally, minimum part-load is modeled with a binary variable $z_i^{\text{on}}$ to ensure the output power $Q_i^{\text{out}}$ is either zero or between minimum part-load $Q_i^{\text{min}}$ and maximum power $Q_i^{\text{max}}$ (\ref{eq:DR_MIMO_P1_part_load}).
Finally, the energy balances state that the demands of the flexible process $Q_{\text{dem},e}^{\text{process}}$ and the demands of other inflexible processes $Q_{\text{dem},e}^{\text{inflexible}}$ have to be satisfied by the output $Q_{i}^{\text{out}}$ of energy conversion units supplying energy $e$, collected in set $\mathbb{C}_e^{\text{sup}}$, and power from the grid $\Delta P_e$ (\ref{eq:DR_MIMO_P1_energy balance}).

In our previous publication \citep{Baader.2021}, we only derived dynamic ramping constraints (\ref{eq:DR_MIMO_P1_dynamic_ramping}) for processes with a single input. We, therefore, only had to constrain a single time derivative of the production rate $\rho$ with order $\delta$, which we defined as the ramping degree of freedom $\nu=\rho^{(\delta)}$. 
In the present contribution, we consider multiple inputs $\vect{u}$ and can potentially have constraints on all considered derivatives of the production rate $\rho^{(\gamma)}$ with $\gamma = 1, ..., \delta$ and the integer $\delta$ being the order of the highest time derivative that is constrained by input bounds.
For instance, the bounds of input $u_1$ could directly limit the first derivative of the production rate, $\dot{\rho}$, whereas the bounds on input $u_2$ could limit the second derivative of the production rate, $\rho^{(2)}$, directly and only influence the first derivative, $\dot{\rho}$, through the integration.
The generalized dynamic ramping constraints (DRCs) developed in this publication therefore read:
\begin{gather}
    \label{eq:DR_MIMO_DRCa}
    \rho^{(\delta)} = \nu \tag{DRCa}
    \\ 
    \label{eq:DR_MIMO_DRCb}
        \dot{\rho}^{\text{min}}(\rho) \leq \dot{\rho} \leq \dot{\rho}^{\text{max}}(\rho) \tag{DRCb}
        \\\nonumber \vdots\\ \tag{DRCc}
        \left(\rho^{(\gamma)} \right) ^{\text{min}}\left(\rho,\dot{\rho},...,\rho^{(\gamma-1)}\right) \leq \rho^{(\gamma)} \leq \left(\rho^{(\gamma)} \right) ^{\text{max}}\left(\rho,\dot{\rho},...,\rho^{(\gamma-1)} \right)
        \\\nonumber \vdots
        \\ 
        \label{eq:DR_MIMO_DRCd} \tag{DRCd}
        \nu^{\text{min}}\left(\rho,\dot{\rho},...,\rho^{(\delta-1)}\right) \leq \nu \leq \nu^{\text{max}}\left(\rho,\dot{\rho},...,\rho^{(\delta-1)} \right)
\end{gather}
Still, only the highest considered time derivative $\rho^{\delta}$ is a degree of freedom and thus defined as the ramping degree of freedom $\nu$.

The limits $\left(\rho^{(\gamma)} \right) ^{\text{min}}$, $\left(\rho^{(\gamma)} \right)^{\text{max}}$ as well as $\nu^{\text{min}}$, $\nu^{\text{max}}$ are in general nonlinear functions that we derive by coordinate transformation (Section~\ref{sec:DRMIMO_method}).
To incorporate the dynamic ramping constraint into an MILP formulation, the true nonlinear limits are approximated by piecewise-affine (PWA) functions for all considered derivatives $\gamma=1,...,\delta$ because PWA functions allow us to formulate a mixed-integer linear scheduling problem.
These piecewise-affine functions need to be conservative such that they prohibit constraint violation with respect to the true nonlinear limits to guarantee that the chosen trajectory for $\rho$ satisfies all bounds on inputs and states. 
Accordingly, the approximations of $\left(\rho^{(\gamma)} \right) ^{\text{min}}$ and $\nu^{\text{min}}$ must always by greater than or equal to the true nonlinear limits and the approximations of $\left(\rho^{(\gamma)} \right) ^{\text{max}}$ and $\nu^{\text{max}}$ must always be lower than or equal to the real nonlinear limits.
Choosing the quality of the approximations allows us to explicitly balance the achievable flexibility range against the computational burden.

The process energy demand $Q_{\text{dem},e}$ (cf. Equation~(\ref{eq:DR_MIMO_P1_process_energy_dem})) for an energy form $e$ is modeled as a function of the production rate and its time derivatives:
    \begin{align}
        \label{eq:DR_MIMO_energy_demand_model}
        Q_{\text{dem},e}^{\text{process}}\left(\rho,\dot{\rho},...,\rho^{(\delta-1)},\nu\right) 
    \end{align}
Similar to the DRC, a piecewise-affine function is chosen for $Q_{\text{dem},e}^{\text{process}}$ to achieve an MILP formulation.

The problem formulation (P) with DRCs captures a larger share of the dynamic flexibility range of processes than static ramping constraints can do while still allowing to formulate a mixed-integer linear program. 
However, to choose suitable piecewise-affine ramping limits, the true nonlinear limits of the process need to be derived or approximated.
In the following section, we show how these limits can be derived rigorously for MIMO processes that are differentially flat.

\section{Deriving dynamic ramping constraints}
\label{sec:DRMIMO_method}
In Section~\ref{subsec:DRMIMO_assumptions}, necessary assumptions are stated  and then our approach is presented in Section~\ref{subsec:DRMIMO_approach}.
\subsection{Assumptions}
\label{subsec:DRMIMO_assumptions}
    \begin{enumerate}
        \item The process degrees of freedom are divided into a control input vector $\vect{u}$ and the variable production rate $\rho$. The process model is a system of ordinary differential equations given by:
        \begin{align}
            \label{eq:DR_MIMO_process_model}
            \dot{\vect{x}} = \vect{f} (\vect{x},\vect{u},\rho) 
        \end{align}
        with state vector $\vect{x} \in \mathbb{R}^n$, and a nonlinear right-hand side function $\vect{f} (\vect{x},\vect{u},\rho)$. We assume that there are no further inputs or disturbances to the process.

        \item The control input vector $\vect{u}$, and the production rate $\rho$ are bounded by minimum and maximum values $\vect{u}^{\text{min}}$, $\vect{u}^{\text{max}}$, and $\rho^{\text{min}}$, $\rho^{\text{max}}$, respectively. Similarly, the states $\vect{x}$ have to be maintained within bounds $\vect{x}^{\text{min}}$, $\vect{x}^{\text{max}}$.
        \item The process (\ref{eq:DR_MIMO_process_model}) is flat. That is, the process has one or multiple flat output vectors $\boldsymbol{\xi}$. An output vector $\boldsymbol{\xi}$ is flat if it satisfies three conditions \citep{Fliess.1995,Rothfu.1997,R.Rothfu.1996}:
            \begin{enumerate}
                \item The flat output vector can be given as a function $\boldsymbol{\phi}$ of states $\vect{x}$, inputs $\vect{u}$, production rate $\rho$, and time derivatives of $\vect{u}$ and $\rho$:
                \begin{align}
                    \label{eq:DR_MIMO_cond_4a}
                    \tag{5a}
                    \boldsymbol{\xi} = \boldsymbol{\phi} \left(\vect{x}, \vect{u}, \dot{\vect{u}}, ..., \vect{u}^{(\alpha)}, \rho, \dot{\rho}, ..., \rho^{(\kappa)}\right)
                \end{align}
                with finite integer numbers $\alpha$, $\kappa$. The function $\boldsymbol{\phi}$ can be seen as a transformation from the original state and input space to the flat output space. Often, it is possible to choose flat outputs that have a physical meaning, e.g., the conversion of a reactor, and that are a function of the states $\vect{x}$ only \citep{Adamy.2014}.
                \item A backtransformation from the flat output and its derivatives to the original states $\vect{x}$ and inputs $\vect{u}$ can be found. Accordingly, the system states $\vect{x}$ and inputs $\vect{u}$ can be given as functions $\boldsymbol{\psi}_1$, $\boldsymbol{\psi}_2$ of the flat outputs $\boldsymbol{\xi}$, the production rate $\rho$, and a number of time derivatives of $\boldsymbol{\xi}$ and $\rho$:
                \begin{align}
                    \tag{5b1}
                    \label{eq:DR_MIMO_cond2}
                    &\vect{x} = \boldsymbol{\psi}_1 \left( \boldsymbol{\xi}, \dot{\boldsymbol{\xi}}, ..., \boldsymbol{\xi}^{(\beta-1)}, \rho, \dot{\rho}, ..., \rho^{(\zeta-1)}\right) \\ 
                    \tag{5b2}
                    &\vect{u} = \boldsymbol{\psi}_2 \left( \boldsymbol{\xi}, \dot{\boldsymbol{\xi}}, ..., \boldsymbol{\xi}^{(\beta)}, \rho, \dot{\rho}, ..., \rho^{(\zeta)}\right)
                \end{align}
                with finite integer numbers $\beta$, $\zeta$. 
                \item The components of $\boldsymbol{\xi}$ are differentially independent \citep{Rothfu.1997}. Consequently, they do not fulfill any differential equation:
                \begin{align}
                \tag{5c}
                    \label{eq:DR_MIMO_cond3}
                    \boldsymbol{\mu}(\boldsymbol{\xi}, \dot{\boldsymbol{\xi}}, ..., \boldsymbol{\xi}^{(\beta)}) = \boldsymbol{0}
                \end{align}
                Condition (\ref{eq:DR_MIMO_cond3}) is satisfied if condition (5b) is satisfied, dim($\boldsymbol{\xi}$) = dim($\vect{u}$) = $m$, and rank$\left(\frac{\partial\vect{f} (\vect{x},\vect{u},\rho)}{\partial \vect{u}}\right)= m$, where $m$ is the number of inputs \citep{Rothfu.1997}.
            \end{enumerate}
            \addtocounter{equation}{1}
            Note: To check conditions (\ref{eq:DR_MIMO_cond_4a}) - (\ref{eq:DR_MIMO_cond3}), a candidate for a flat output vector $\boldsymbol{\xi}$ is needed. We assume that such a candidate for a flat output vector $\boldsymbol{\xi}$ can be identified based on engineering intuition. 
        \item The trajectory of the production rate $\rho$ is determined by the scheduling optimization. Subsequently, the control input vector $\vect{u}\in\mathbb{R}^m$ is calculated by an underlying process control.
    \end{enumerate}

\begin{figure*}[h]
    \centering
 \includegraphics{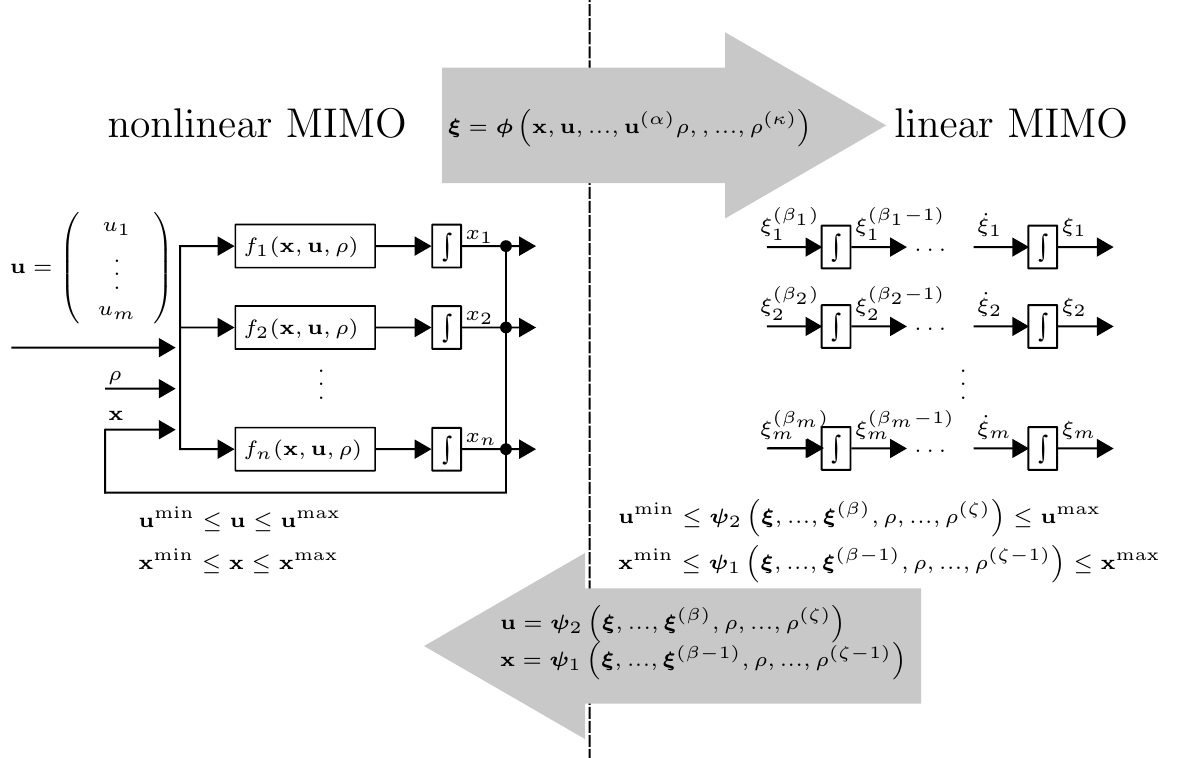}
  \caption{Visualization of flatness-based coordinate transformation showing original nonlinear MIMO process model as introduced in assumption 1 (left), transformations $\boldsymbol{\phi}, \boldsymbol{\psi}_1, \boldsymbol{\psi}_2$ as introduced in assumption 3 (gray arrows), and linear MIMO process model in transformed coordinate space as introduced in assumption 3 (right).}
   \label{fig:DRMIMO_concept_flatness}
\end{figure*}

The flatness-based coordinate transformation is visualized in Figure \ref{fig:DRMIMO_concept_flatness}.
A nonlinear model with linear constraints is transformed into a linear model with nonlinear constraints.
The number of control degrees of freedom is maintained: While, in the original model, the $m$ inputs $u_k$ ($k=1,...,m$) are the degrees of freedom, in the linear model, the $m$ highest time derivatives of the output components $\xi_k^{(\beta_k)}$ ($k=1,...,m$) are the degrees of freedom \citep{Adamy.2014,Fliess.1995}.
The number of time derivatives $\beta_k$ can deviate between the flat output components $\xi_k$ \citep{Adamy.2014,Fliess.1995}.
In the linear model, the transformed state vector $\boldsymbol{\Xi}$ is formed by the outputs and their time derivatives (except for the highest time derivative):
\begin{align}
    \label{eq:DRMIMO_Xi}
    \boldsymbol{\Xi} = \left(\xi_1, ...,\xi_1^{(\beta_1-1)},\xi_2, ...,\xi_2^{(\beta_2-1)}, ..., \xi_m, ... ,\xi_m^{(\beta_m-1)}  \right)^T
\end{align}
The dimension of the transformed state vector $\boldsymbol{\Xi}$ is greater or equal to the dimension of the original state vector $\vect{x}$ \citep{Adamy.2014,Fliess.1995}. 
Consequently, the original nonlinear process model is converted to a linear model consisting of $m$ integrator chains.
In the linear model, every flat output $\xi_k$ can be varied with a $\beta_k$-th order dynamic independently of the other flat outputs.
Note that flatness is a sufficient condition for controllability of a nonlinear process \citep{Adamy.2014}.

To ease notation in the following, we introduce the ramping state vector 
\begin{align}
   \boldsymbol{\varphi} = \left(\rho, \dot{\rho}, \dots ,\rho^{(\delta-2)}, \rho^{(\delta-1)}\right)^T
\end{align}
and its time derivative 
\begin{align}
    \dot{\boldsymbol{\varphi}} = \left(\dot{\rho}, \rho^{(2)}, \dots,\rho^{(\delta-1)}, \nu \ \right)^T.
\end{align}

\subsection{Approach}
\label{subsec:DRMIMO_approach}

\begin{figure*}[h]
\centering 
 \includegraphics{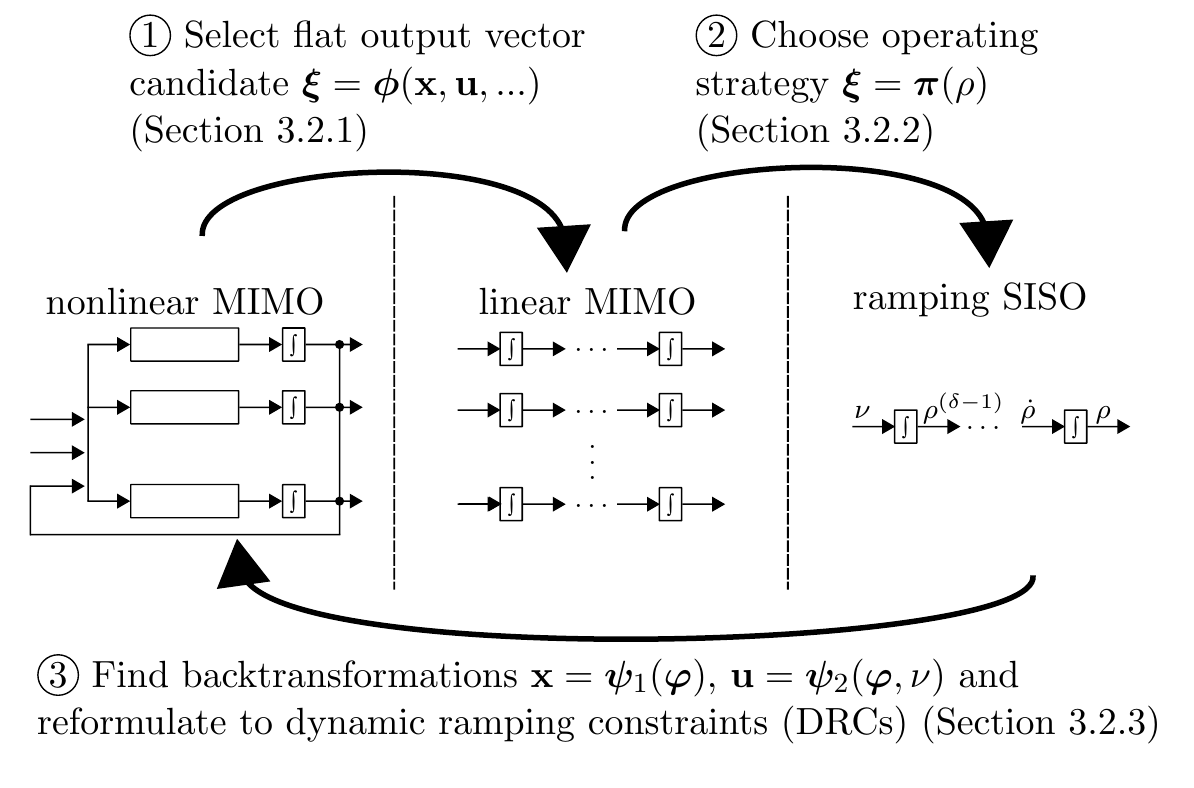}
  \caption{Overview of steps to derive dynamic ramping constraints performed in the respective sections. For the nonlinear MIMO model and the linear MIMO model, all variable symbols are omitted for clarity as they are identical to Figure 2. The variables of the ramping SISO model are the production rate $\rho$, its derivatives $\dot{\rho},...,\rho^{(\delta-1)}$, and the ramping degree of freedom $\nu$.}
   \label{fig:DRMIMO_steps}
\end{figure*}

The steps to derive dynamic ramping constraints further discussed in the following subsections are summarized in Figure~\ref{fig:DRMIMO_steps}: 
First, a candidate for a flat output vector is selected based on two necessary flatness conditions (Section~\ref{sec:DRMIMO_necess_cond}).
Assuming this candidate is, in fact, a flat output, the original nonlinear MIMO process model is transformed into a linear  MIMO model (center of Figure~\ref{fig:DRMIMO_steps}).
In this linear model, the $m$ components of the flat output vector $\xi_k$ are decoupled such that the outputs $\xi_k$ can be controlled independently of each other by manipulating the degrees of freedom $\xi_k^{(\beta_k)}$ (Figure \ref{fig:DRMIMO_concept_flatness}). 
However, a SISO model is needed for the dynamic ramping constraints where the production rate $\rho$ is controlled by manipulating the ramping degree of freedom $\nu$. 
Thus, second, the flat output components $\xi_k$ are coupled by choosing an operating strategy $\boldsymbol{\pi}(\rho)$ that defines every $\xi_k$ as function of the production rate $\rho$ (Section~\ref{sec:DRMIMO_op_strategy}).
This coupling reduces the number of degrees of freedom from $m$ to one, leading to a significant model order reduction.
Third, backtransformations are found from the ramping SISO model to the original nonlinear MIMO model, and thus flatness is proven (Section~\ref{sec:DRMIMO_reformulation}).

Based on the backtransformations and the ramping SISO model, the true nonlinear ramping limits are derived from the bounds on inputs $\vect{u}$ and states $\vect{x}$. After approximating the true nonlinear limits with piecewise-affine functions, the dynamic ramping constraints are complete and the problem (P) can be solved as MILP.

\subsubsection{Selection of flat output candidate and necessary flatness conditions}
\label{sec:DRMIMO_necess_cond}

First, we propose to apply a necessary condition for flatness from literature based on a graph representation of the process \citep{Schulze.2020} to identify a candidate for the flat output vector $\boldsymbol{\xi}$ as  $\boldsymbol{\phi} \left(\vect{x}, \vect{u}, \dot{\vect{u}}, ..., \vect{u}^{(\alpha)}, \rho, \dot{\rho}, ..., \rho^{(\kappa)}\right)$ (compare to condition 5a).
As the second necessary condition for flatness, we propose setting up an equation system that implicitly defines the backtransformation $\boldsymbol{\psi}$ (compare to condition 5b) and check if this backtransformation equation system is structurally solvable.

For the graph representation, the process model is represented as a directed graph in which all states $x_i$ and inputs $u_k$ are represented as vertices $v_{x_i}$ and $v_{u_k}$, respectively \citep{Schulze.2020}.
If $\frac{\partial f_i(\vect{x},\vect{u})}{\partial x_j}$ is non-zero, there is an edge from vertex $v_{x_j}$ to vertex $v_{x_i}$ \citep{Schulze.2020}.
In other words: If a state $x_j$ acts on the derivative of another state $x_i$, an edge is drawn from $x_j$ to $x_i$ .
Similarly, an edge from input $u_k$ to state $x_i$ exists, if $\frac{\partial f_i(\vect{x},\vect{u})}{\partial u_k}$ is non-zero.
The necessary condition for a flat output vector $\boldsymbol{\xi}$ is that it must be possible to match the $m$ output components $\xi_i$ to the input components $u_i$ such that there are $m$ input-output pairs with pairwise disjoint paths through the graph that cover all process states \citep{Schulze.2020}. 
An illustrative example for this necessary condition is given in the Supplementary Information (SI).

As a starting point, the states which are typically controlled can be tested as flat output candidates.
Typical states to be controlled are outlet stream compositions, final ef\/f\/luent temperature, or the process hold-up \citep{Jogwar.2009}.

Once a flat output vector candidate is identified based on the graph representation as $\boldsymbol{\xi} = \boldsymbol{\phi} (\vect{x}, \vect{u}, \dot{\vect{u}}, ...,$ $\vect{u}^{(\alpha)}, \rho, \dot{\rho}, ..., \rho^{(\kappa)})$, condition 4a is fulfilled.
Further, the rank criterion needed for condition 4c can be checked easily.
As a final step to show flatness, existence of the backtransformations $\boldsymbol{\psi}_1, \boldsymbol{\psi}_2$ that give states $\vect{x}$ and inputs $\vect{u}$ as functions of the flat output $\boldsymbol{\xi}$ and its time derivatives $\dot{\boldsymbol{\xi}}, ..., \boldsymbol{\xi}^{(\beta)}$ needs to be shown (condition 4b).
To obtain these backtransformations, a nonlinear system of equations needs to be developed with $\boldsymbol{\xi}$ and its derivatives $\dot{\boldsymbol{\xi}}, ..., \boldsymbol{\xi}^{(\beta)}$ on the right-hand side and left-hand side functions including the wanted quantities: states $\vect{x}$, inputs $\vect{u}$, and potentially derivatives of inputs.
To set up the backtransformation, we make use of the fact that $\boldsymbol{\xi}$ is given by $\boldsymbol{\phi} \left(\vect{x}, \vect{u}, \dot{\vect{u}}, ..., \vect{u}^{(\alpha)}, \rho, \dot{\rho}, ..., \rho^{(\kappa)}\right)$ and the derivatives of $\boldsymbol{\xi}$ can be given as functions of $\vect{x}$, inputs $\vect{u}$, and  derivatives of inputs by means of the total differential.
The equation system for the backtransformation needs to be square such that there are as many equations as there are unknown states, inputs, and derivatives of inputs.
We propose to check if this nonlinear equation system is structurally solvable by conducting an analysis similar to the structural index analysis for differential-algebraic equation systems \citep{J.Unger.1995} commonly used in process systems engineering.
An example is given in the SI.

\subsubsection{Operating strategy}
\label{sec:DRMIMO_op_strategy}
Under the assumption that the flat output candidate identified as discussed in Section~\ref{sec:DRMIMO_necess_cond} is, in fact, a flat output, the nonlinear model can be transformed into a linear model.  
This linear model is a MIMO model with $m$ degrees of freedom $\xi_i^{(\beta_i)}$ and $m$ outputs $\xi_i$ (Figure \ref{fig:DRMIMO_concept_flatness}).
As discussed above, the linear MIMO is transferred to a SISO ramping model by coupling the flat output components.  
To this end, we insert an operating strategy that gives the value of every flat output as a function of the production rate $\rho$. 
This operating strategy differentiates between two possible types of flat output components: 
First, flat output components might have specifications that should be maintained constant, such as outlet stream compositions, final ef\/f\/luent temperature, and the process hold-up \citep{Jogwar.2009}.
Accordingly, the operating strategy is to hold such an output  $\xi_k$ constant at its nominal value such that $\xi_k=\xi_k^{\text{nom}}$ holds and thus all time derivatives are zero.

Second, flat output components may be unspecified.
For instance, in our case study, one flat output component is a concentration for which no specifications are given.
As every flat output component $\xi_k$ corresponds to one control input, i.e., degree of freedom, $u_k$, if specifications are given for $l$ outputs, and $l$ is smaller than $m$, $m-l$ flat output components remain as degrees of freedom in steady-state.
Thus, the $m-l$ free flat output component $\xi_k$  can, in principle, be chosen to have any value in steady-state as long as no variable bounds are violated.
To have the optimal steady-state operating points, we use a steady-state optimization to determine the optimal values $\xi_k$ in advance as a function of the production rate $\rho$ such that $\xi_k = \pi_k(\rho)$.
For instance, the objective can be to find the steady-state operating points that minimize energy consumption. 
In our case study, we choose the flat output component such that the overall heat demand is minimal for steady-state points (Section~\ref{sec:DRMIMO_case_study}).

The operating strategy $\pi_k(\rho)$ can be any nonlinear function.
The only requirement is that the function $\pi_k(\rho)$ must be differentiable with respect to $\rho$ sufficiently often so that all derivatives of $\xi_k$ which are part of the backtransformation discussed in the previous section are defined by the total differential, e.g., $\dot{\xi}_k = \frac{\partial\pi_k(\rho)}{\partial \rho}\dot{\rho}$, $\xi_k^{(2)} = \frac{\partial\pi_k(\rho)}{\partial \rho}\rho^{(2)} + \frac{\partial^2\pi_k(\rho)}{\partial\rho^2 }\dot{\rho}^2$.

When all output components with specifications are maintained at their nominal values and a function $\pi_k(\rho)$ is chosen for all other output components, the operating strategy can be summarized as $\boldsymbol{\xi}=\boldsymbol{\pi}(\rho)$.
Consequently, the flat output vector $\boldsymbol{\xi}$ and all relevant time derivatives $\boldsymbol{\xi}, \dot{\boldsymbol{\xi}}, ..., \boldsymbol{\xi}^{(\beta)}$ (compare to equation (5b)) are defined as function of the production rate $\rho$ and a number of its time derivatives.
The highest time derivative $\rho^{(\delta)}$ that occurs defines the order of the dynamic ramping constraint (DRC) and the ramping degree of freedom $\nu = \rho^{(\delta)}$.
Consequently, the backtransformations $\boldsymbol{\psi}_1, \boldsymbol{\psi}_2$ discussed in the following only depend on the ramping state vector $\boldsymbol{\varphi}$ and the ramping degree of freedom $\nu$.

\subsubsection{Backtransformation and reformulation to dynamic ramping constraints}
\label{sec:DRMIMO_reformulation}
First, the flat output vector $\boldsymbol{\xi}$ and its time derivatives $\boldsymbol{\xi}, \dot{\boldsymbol{\xi}}, ..., \boldsymbol{\xi}^{(\beta)}$ are replaced in the nonlinear equation system for the backtransformation derived in Section~\ref{sec:DRMIMO_necess_cond}. 
While $\boldsymbol{\xi}$ is replaced by $\boldsymbol{\pi}(\rho)$, the derivatives are replaced by building the total differential of $\boldsymbol{\pi}(\rho)$.
Next, we solve the equation system to get $\vect{x}=\boldsymbol{\psi}_1(\boldsymbol{\varphi})$ and $\vect{u}=\boldsymbol{\psi}_2(\boldsymbol{\varphi},\nu)$.
It is favorable to solve the equation system analytically.
Still, it is not necessary to derive the functions $\boldsymbol{\psi}_1, \boldsymbol{\psi}_2$ analytically as they can also be evaluated numerically as long as their solution is unique.
In this paper, we use the computer algebra package SymPy \citep{meurer2017sympy} to obtain analytic functions $\boldsymbol{\psi}_1, \boldsymbol{\psi}_2$.
Note that it might not be possible to solve the system of equations because the graphical and structural criteria proposed in Section~\ref{sec:DRMIMO_necess_cond} are only necessary flatness conditions. In case the nonlinear system of equations cannot be solved, one can test other flat output candidates.
If the functions $\boldsymbol{\psi}_1, \boldsymbol{\psi}_2$ are found, condition 3b is fulfilled, and flatness is shown at least locally on the subspace defined by the operating strategy.

This flat system constitutes one integrator chain with the degree of freedom $\nu$ that can be chosen arbitrarily at any point in time.
Thus, mathematically, the production rate can be changed infinitely fast by choosing infinitely high values for $\nu$.
However, the real process control inputs $\vect{u}$ are bounded by maximum and minimum values (assumption 3), and therefore, dynamic ramping constraints (DRCs) are needed to ensure that the real process inputs are maintained within bounds.

To get the dynamic ramping constraints (DRCs), we consider the input bounds row by row.
For the $k$-th row, we get
    \begin{align}
        \label{eq:DR_MIMO_u_k}
        u_k^{\text{min}} \leq \psi_{2,k}(\rho,\dot{\rho}, ...,\rho^{(\gamma)}) \leq u_k^{\text{max}} ,
    \end{align}
    where $\rho^{\gamma}$ is the highest time derivative needed to compute $u_k$.
Equation~(\ref{eq:DR_MIMO_u_k}) implicitly limits $\rho^{\gamma}$ for given $\rho,\dot{\rho},..,\rho^{(\gamma-1)}$. 
These limits can be given as an explicit function if it is possible to symbolically invert $\psi_{2,k}$ for the derivative ${\rho^{(\gamma)}}_k$ and thus derive an analytic function $\theta$ of the production rate $\rho$, derivatives $\dot{\rho}, ...,\rho^{\gamma-1}$, and the input $u_k$:
    \begin{align}
        {\rho^{(\gamma)}}_k = \theta\left(\rho,\dot{\rho}, ...,\rho^{\gamma-1}, u_k\right)
    \end{align}
Inserting $u_k^{\text{min}}$, $u_k^{\text{max}}$ into $\theta$ gives ${\rho^{(\gamma)}}_k^{\text{min}}$, ${\rho^{(\gamma)}}_k^{\text{max}}$.
Alternatively, the limits can be derived numerically for given $\rho,\dot{\rho}, ...,\rho^{\gamma-1}$ by sampling different values for ${\rho^{(\gamma)}}_k$ and evaluating if the resulting $u_k$ is within the allowed range.

Note: Here, we discuss the case of one allowed region for ${\rho^{(\gamma)}}_k$ given by limits ${\rho^{(\gamma)}}_k^{\text{min}}$, ${\rho^{(\gamma)}}_k^{\text{max}}$. 
In principle, there could be two (or even more) non-connected allowed regions given by limits  ${\rho^{(\gamma)}}_k^{\text{min,1}}$, ${\rho^{(\gamma)}}_k^{\text{max,1}}$ and ${\rho^{(\gamma)}}_k^{\text{min,2}}$, ${\rho^{(\gamma)}}_k^{\text{max,2}}$ with ${\rho^{(\gamma)}}_k^{\text{max,1}} < {\rho^{(\gamma)}}_k^{\text{min,2}}$, as $\psi_{2,k}$ is a general nonlinear function. In that case, one could either restrict the operating range to one of the two regions or introduce a binary variable that indicates which region is active. 
An illustrative example that can result in two allowed regions is given in the SI.

Bounds on a state $x_k$ give an equation which is has the same structure as
Equation~(\ref{eq:DR_MIMO_u_k}).
Thus, constraints on states can be treated in the same way as constraints on inputs.

\begin{figure*}[h]
    \centering
 \includegraphics[width=8cm]{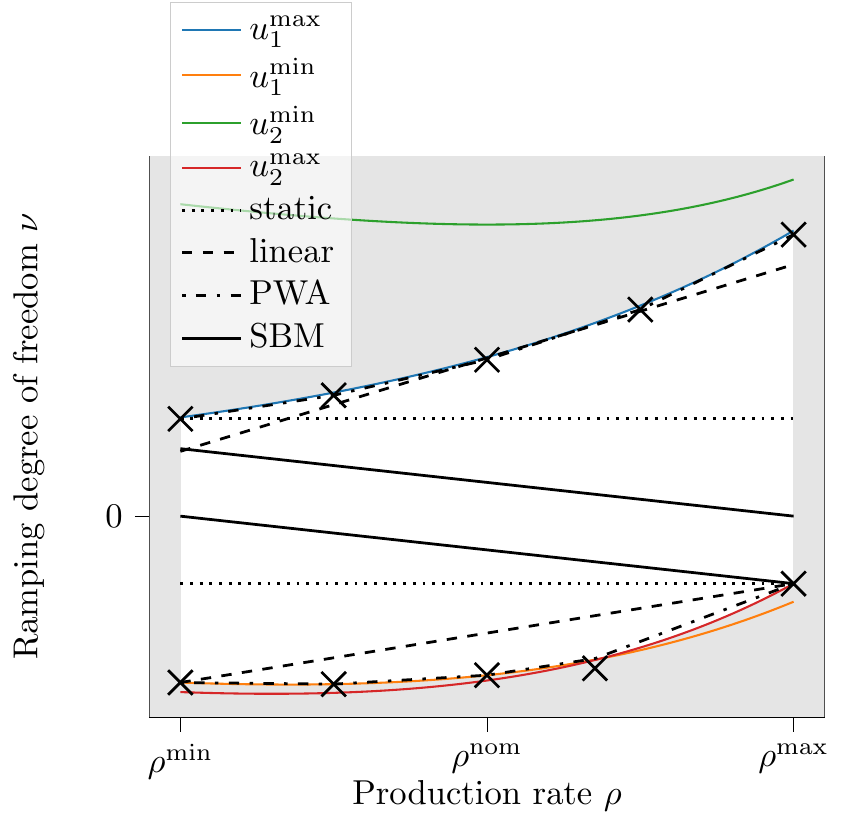}
  \caption{Constraints for ramping degree of freedom $\nu$ as function of production rate $\rho$ for an illustrative case with first-order ramping constraints and two limiting inputs $u_1$, $u_2$. For first-order ramping constraints, the ramping state vector $\boldsymbol{\varphi}$ is of dimension one and equal to the production rate $\rho$. Consequently, the limits on the ramping degree of freedom $\nu$ only depend on $\rho$. The true nonlinear limits caused by the minimum and maximum values of the two inputs $u_1$, $u_2$ are compared to static limits (dotted), linear limits (dashed) and piecewise-affine (PWA) limits (dashed-dotted). Moreover, a linear scale-bridging model (SBM) is visualized (compare to discussion in Section~3.3)}
   \label{fig:DRMIMOSI_comp_bounds_MIMO}
\end{figure*}

Finally, the nonlinear ramping limits need to be approximated by piecewise affine (PWA) functions.
This approximation allows to explicitly balance the quality of the dynamic ramping constraints against the computational burden.
In contrast to previous work \citep{Baader.2021}, there can more be several inputs and states limiting the same derivative of the production rate $\rho$ for the MIMO case. 
Figure~\ref{fig:DRMIMOSI_comp_bounds_MIMO} shows an illustrative case with first-order ramping and two limiting inputs $u_1$, $u_2$.
The upper ramping limit is only determined by the maximum input $u_1^{\text{max}}$ as the upper limit in the flat system resulting from $u_2^{\text{min}}$ is always above the limit from $u_1^{\text{max}}$.
In contrast, the lower ramping limit is given by an intersection between the lower limits from $u_1^{\text{min}}$ and $u_2^{\text{max}}$, respectively. 
If static ramping limits are chosen, a large amount of the feasible region needs to be cut off for the illustrative case in Figure \ref{fig:DRMIMOSI_comp_bounds_MIMO} (horizontal dotted lines).
In contrast, linear (dashed lines in Figure \ref{fig:DRMIMOSI_comp_bounds_MIMO}) and piecewise affine (dashed-dotted lines in Figure \ref{fig:DRMIMOSI_comp_bounds_MIMO}) limits allow to come closer to the true nonlinear limits and thus realize a larger flexibility range.
For the lower ramping limit, piecewise affine limits can be realized without the addition of binary variables as the feasible region is convex.
However, for the upper limit, the feasible region is non-convex and binary variables are needed, making the optimization more computationally challenging.

As, in the general case, the limit functions $\left(\rho^{(\gamma)} \right) ^{\text{min}}$, $\left(\rho^{(\gamma)} \right)^{\text{max}}$ are multivariate functions, multivariate regression methods, e.g., hinging hyperplanes \citep{L.Breiman.1993,A.A.Adeniran.2017,Kamper.2021b}, convex region surrogates \citep{Zhang.2016b,Schweidtmann.2021}, or artificial neural networks with ReLU activation functions \citep{BjarneGrimstad.2019,reluMIP.2021}, can be used to find piecewise-affine approximations.

\subsection{Comparison to other approaches}
\label{sec:comp}
Finally, we compare our approach to two other relevant approaches that integrate scheduling and control by considering a simplified version of the process dynamics in scheduling: scale-bridging models \citep{Du.2015} and data-driven closed-loop models \citep{Kelley.2018}. 
The difference between the two alternative approaches is that scale-bridging models explicitly adapt the underlying control to linearize the closed-loop response \citep{Du.2015} while data-driven closed-loop models identify the response of the process to a change of the set-point $\rho_{SP}$ from data \citep{Dias.2016,Diangelakis.2017,Burnak.2018,Pattison.2016}.
Thus, data-driven closed-loop models, in general, identify a nonlinear closed-loop response, and scale-bridging models rely on a linearization performed by the underlying control.
This linearization can be achieved by exact input-output linearization control \citep{Du.2015}, scheduling-oriented model-predictive control \citep{Baldea.2015}, or by a combination of a set-point filter and tracking control \citep{baader2022simultaneous}.
For first-order dynamics, the linearized closed-loop response reads
\begin{align}
    \label{eq:SBM}
    \rho + \tau \dot{\rho} = \rho_{\text{SP}}
\end{align}
In Equation~(\ref{eq:SBM}), $\tau$ is a tunable time constant and $\rho_{\text{SP}}$ is the set-point given to the underlying controller. That is, instead of the ramping degree of freedom $\nu$, the set-point $\rho_{\text{SP}}$ is the degree of freedom for the scheduling optimization. 
To compare the scale-bridging model to our dynamic ramping constraints, we rearrange Equation~(\ref{eq:SBM}) to calculate the maximum possible rate of change $\dot{\rho}^{\text{max}}$ as a function of the maximum set-point $\rho_{SP}^{\text{max}}$ and the minimum rate of change $\dot{\rho}^{\text{min}}$ as a function of the minimum set-point $\rho_{SP}^{\text{min}}$:
\begin{align}
    \label{eq:SBM_rho_max}
    \dot{\rho}^{\text{max}} = \frac{1}{\tau}(\rho_{SP}^{\text{max}} - \rho) \\
    \label{eq:SBM_rho_min}
    \dot{\rho}^{\text{min}} = \frac{1}{\tau}(\rho_{SP}^{\text{min}} - \rho)
\end{align}
The resulting limits are visualized in Figure~\ref{fig:DRMIMOSI_comp_bounds_MIMO} for the natural choice that the maximum set-point $\rho_{SP}^{\text{max}}$ equals the maximum production rate $\rho^{\text{max}}$ and the minimum set-point $\rho_{SP}^{\text{min}}$ equals the maximum production rate $\rho^{\text{min}}$. The time constant $\tau$ is chosen such that the rate of change $\dot{\rho}$ demanded by the scale-bridging model is always within limits.
The flexibility range for the scale-bridging model has a parallelogram shape dictated by Equations~(\ref{eq:SBM_rho_max}) and (\ref{eq:SBM_rho_min}).
Our dynamic ramping constraints can always match this parallelogram shape by linear ramping limits $\nu^{\text{min}}$, $\nu^{\text{max}}$.
Additionally, dynamic ramping constraints also allow choosing linear limits with another shape or even piecewise affine limits.
Thus dynamic ramping constraints can consistently perform at least as well as linear scale-bridging models.
Moreover, our analysis could be used to rigorously choose the time constant $\tau$ of the scale-bridging model.

As stated above, data-driven closed-loop models identify a nonlinear closed-loop response.
Still, if piecewise affine functions approximate nonlinearities, a MILP formulation is possible, see, e.g., \cite{Kelley.2018} who derive an MILP formulation for Hammerstein-Wiener models.
Thus, data-driven closed-loop models are an alternative to our dynamic ramping constraints; in particular, they are still applicable if no mechanistic process model is available or the process is not flat. 
If a flat mechanistic process model is available, there are two conceptual differences between dynamic ramping constraints and data-driven closed-loop models: 
 First, dynamic ramping constraints have the theoretical advantage that they can be derived rigorously such that the feasibility of the optimized trajectory is guaranteed. 
Second, dynamic ramping constraints and data-driven closed-loop models lead to conceptually different trade-offs between computational burden and accuracy: While dynamic ramping constraints reduce the computational burden by reducing the accuracy of the ramping limits, data-driven closed-loop models reduce computational burden by reducing the accuracy of the model prediction.
Thus, dynamic ramping constraints compromise on the feasible region, sacrificing flexibility, and data-driven closed-loop models compromise on the prediction of the process outputs, potentially leading to constraint violations or production deviating from the schedule.

\section{Case study: Reactor-separator process with recycle}
\label{sec:DRMIMO_case_study}

\begin{figure*}[h]
\centering
 \includegraphics{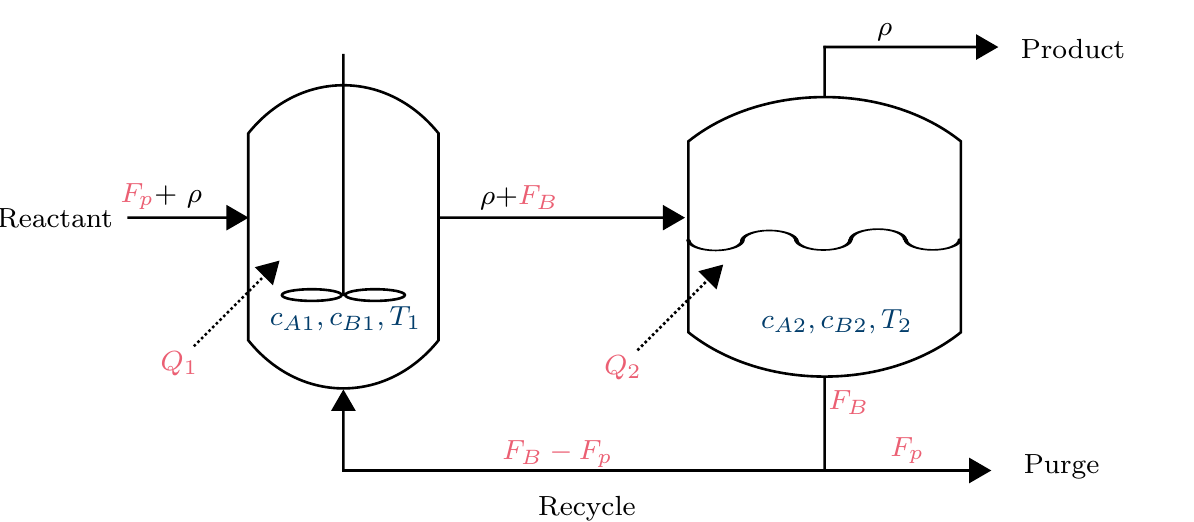}
  \caption{Case study of reactor-separator process with recycle: States \textcolor{JB}{$\vect{x}$} are the concentrations of component $A$ and $B$, \textcolor{JB}{$c_{A}$}, and \textcolor{JB}{$c_{B}$}, respectively, and the temperature \textcolor{JB}{$T$} in the reactor (1) and in the flash (2). 
  Manipulated control inputs \textcolor{JR}{$\vect{u}$} are the bottom stream \textcolor{JR}{$F_B$}, the purge stream \textcolor{JR}{$F_P$}, the heat input to the reactor \textcolor{JR}{$Q_1$}, and the heat input to the flash \textcolor{JR}{$Q_2$}.
  The scheduling degree of freedom is the production rate $\rho$.
  All other material flow rates are given as functions of \textcolor{JR}{$F_B$}, \textcolor{JR}{$F_p$}, and $\rho$, e.g., the reactant stream is equal to the purge \textcolor{JR}{$F_p$} plus the production rate $\rho$ as no accumulation of material occurs. }
   \label{fig:DRMIMO_case_study}
\end{figure*}

In this case study, we consider a reactor-separator process with recycle consisting of a continuous stirred tank reactor (CSTR) and a flash (Figure \ref{fig:DRMIMO_case_study}).
The production rate $\rho$ can be varied around its nominal value $\rho^{\text{nom}}$ between $\rho^{\text{min}} = 0.8\rho^{\text{nom}}$ and $\rho^{\text{max}} = 1.2\rho^{\text{nom}}$ as long as the nominal production is reached on average over the considered time horizon.
A raw material $A$ reacts to the desired product $B$, which can further react to an undesired product $C$.
The process has 6 differential states: the concentration of $A$, $c_{A1}$, the concentration of $B$, $c_{B1}$, and the temperature $T_1$ in the reactor (1) and the analog quantities $c_{A2}$, $c_{B2}$, $T_2$, in the flash (2). 
Apart from the production rate $\rho$, there are four manipulated variables: the bottom stream $F_B$, the purge stream $F_p$, the heat flow to the reactor $Q_1$, and the heat flow to the flash $Q_2$.
The process equations are modified from the textbook example by \cite{Christofides.2011}~where 2 CSTRs and a flash are considered.
Though, also the original version \citep{Christofides.2011} fulfills our assumptions, we decided to modify the example to reduce the number of states from 9 to 6 to improve readability and clarity. 

The model equations comprise component and energy balances given in the SI. 

\subsection{Selection of flat output candidate and necessary flatness conditions}
\begin{figure}[h]
\centering
\includegraphics{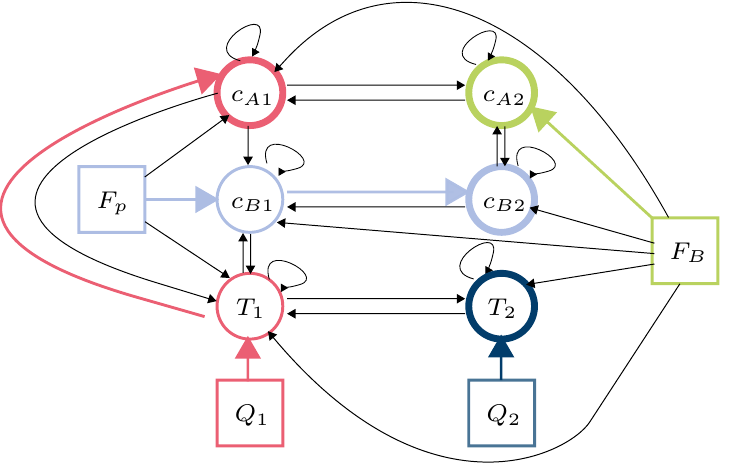}
  \caption{Graph representation for the reactor-separator process with recycle (compare to Figure~\ref{fig:DRMIMO_case_study}). The output $\boldsymbol{\xi}=(c_{A2}, c_{B2}, T_2, c_{A1})^T$ fulfills the necessary condition for a flat output.}
   \label{fig:DRMIMO_graph_case_study}
\end{figure}

A flat output candidate $\boldsymbol{\xi}$ must have 4 components as there are four control inputs.
We first consider the three states of the flash $c_{A2}$, $c_{B2}$, and $T_2$, as they determine the outlet stream, and we assume that specifications for the outlet stream are given.
As fourth output component, we choose the concentration $c_{A1}$ based on the graph representation in Figure \ref{fig:DRMIMO_graph_case_study}.
The four input-output pairs $F_B - c_{A2}$, $F_p - c_{B2}$, $Q_2 - T_2$, $Q_1 - c_{A1}$ fulfill the necessary condition for a flat output (Figure \ref{fig:DRMIMO_graph_case_study}).

By differentiating the components of the output $\boldsymbol{\xi}=\left(c_{A2}, c_{B2}, T_2, c_{A1} \right)$ up to three times, we receive a structurally solvable system of equations (Table~S3 in the SI).

\subsection{Operating strategy}
\label{subsec:DRMIMO_case_study_stratey}
For the operating strategy, we assume that the composition and temperature of the product stream $\rho$ must be maintained constant.
Accordingly, $\xi_1=c_{A2}$, $\xi_2=c_{B2}$, $\xi_3=T_2$ must be maintained at their nominal values $\xi_1^{\text{nom}} = 0.4539$, $\xi_2^{\text{nom}} = 0.4610$, $\xi_3^{\text{nom}}=455 K$. 
Thereby, the considered operating region is already significantly reduced. 

As there are 4 control inputs, we can maintain the first three flat output components at their nominal values and still have one degree of freedom left.
Consequently, an operating strategy for the fourth flat output $\xi_4=c_{A1}$ can be chosen freely.
In a steady-state optimization, we search for the steady-state operating points that minimize the total heating $Q_1 + Q_2$ and fix $\xi_4$ to be a function of the production rate $\xi_4=\pi_4(\rho)$.
Further details on this steady-state optimization are given in the SI.

\subsection{Reformulation to dynamic ramping constraints}
After inserting the operating strategy defined in Section~\ref{subsec:DRMIMO_case_study_stratey}, we solve the nonlinear backtransformation equation system.
To this end, we use the computer algebra package SymPy \citep{meurer2017sympy} to find explicit algebraic expressions for all states and inputs, except for the temperature $T_1$.
The temperature $T_1$ has to be determined numerically because the two reaction terms in the differential Equation~(S9) lead to an equation of the form
\begin{align}
    \label{eq:DR_MIMO_case_study_numeric}
    0= b_0 + b_1e^{-\frac{\textcolor{gray}{E_1}}{\textcolor{gray}{R}T_1}} + b_2 e^{-\frac{\textcolor{gray}{E_2}}{\textcolor{gray}{R}T_1}},
\end{align}
with parameters $b_0$, $b_1$, $b_2$.
The implicit Equation~(\ref{eq:DR_MIMO_case_study_numeric}) has a unique solution for the temperature, as the exponential functions are monotonic. 
The numerical solution of Equation~(\ref{eq:DR_MIMO_case_study_numeric}) is found using the python package SciPy \citep{2020SciPy-NMeth}.

Table \ref{tab:DRMIMO_case_study_dependency_rho} provides the structural dependency of the states and inputs on the production rate $\rho$ and its time derivatives $\dot{\rho}$ and $\nu$ ($=\rho^{(2)}$).
The highest time derivative of the production rate that appears is two.
Therefore, the ramping degree of freedom $\nu$ is equal to $\rho^{(2)}$.
More detailed information about the reformulations is given in the SI.

\begin{table}[bt]
    \centering
    \setlength{\tabcolsep}{3pt}
    \renewcommand{\arraystretch}{1.0}
    \caption{Structural dependency resulting from nonlinear transformation and operating strategy of states and inputs on the production rate $\rho$ and its time derivatives $\dot{\rho}$ and $\nu$ ($=\rho^{(2)}$). States which are held constant are marked with a star (*).}
    \label{tab:DRMIMO_case_study_dependency_rho}
    \begin{tabular}{c|c|c|c|c}
    \multicolumn{1}{c}{}&& $\rho$ & $\dot{\rho}$ & $\nu$ \\
    \hline
    \multirow[t]{6}{*}{states:}  & $c_{A1}$ & x & & \\
      &$c_{B1}$ & x & & \\
     &$T_1$ & x & x & \\
    &$c_{A2}$* &  & & \\
    &$c_{B2}$* &  & & \\
    &$T_2$* &  & & \\
    \hline
    \multirow[t]{4}{*}{inputs:}& $F_B$ & x &  & \\
    &$F_p$ &x   &x  & \\
    &$Q_{1}$ &x  &x &x \\
    &$Q_2$ &x  &x & \\

    \end{tabular}
\end{table}

For the dynamic ramping constraints, it has to be analyzed which state and input bounds limit the time derivatives of the production rate $\rho$.
States and inputs that are either held constant or exclusively depend on the production rate $\rho$ do not need to be checked because it is already known from the steady-state optimization that these variables take feasible values for all considered production rates.
Consequently, the variables $c_{A,1}, c_{B,1}, c_{A,2}, c_{B,2}, T_2, F_B$ do not influence the dynamic ramping constraints (compare to Table~\ref{tab:DRMIMO_case_study_dependency_rho}).

\begin{figure}[h]
    \centering
 \includegraphics{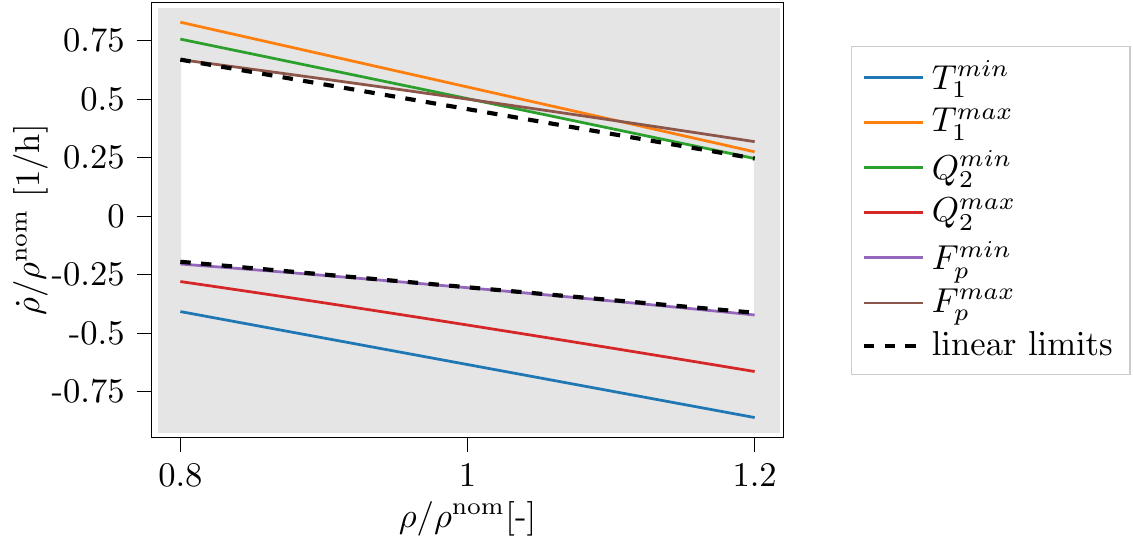}
  \caption{Limits on the scaled time derivative of the production rate $\dot{\rho}$ over the scaled production rate $\rho$ resulting from the minimum (min) and maximum (max) values of the state $T_1$ and the inputs $F_p$, $Q_2$.}
   \label{fig:DRMIMO_bounds_rho_1}
\end{figure}

The variables $T_1, F_p, Q_2$ depend on $\rho$ and $\dot{\rho}$ but not on the second time derivative $\nu$.
Accordingly, the bounds of these variables limit the first time derivative $\dot{\rho}$.
In Figure \ref{fig:DRMIMO_bounds_rho_1}, the limits on $\dot{\rho}$ resulting from variable bounds are shown over the production rate $\rho$.
The calculation of these limits is explained in the SI.
While the lower limit on $\dot{\rho}$  results from the bound $F_p^{\text{min}}$, the upper limit is given by $F_p^{\text{max}}$ for small production rates and by $Q_2^{\text{min}}$ for large production rates.
Graphically, we choose conservative linear functions for the limits $\dot{\rho}^{\text{min}}(\rho)=a_0^{\text{min}} + a_1^{\text{min}}\rho$, $\dot{\rho}^{\text{max}}(\rho)=a_0^{\text{max}} + a_1^{\text{max}}\rho$ with parameters $a_0^{\text{min}}$, $a_1^{\text{min}}$, $a_0^{\text{max}}$, $a_1^{\text{max}}$ (compare to (\ref{eq:DR_MIMO_DRCb})).

The heating input $Q_1$ is the only variable that depends on the ramping degree of freedom $\nu$. 
We choose piecewise affine limits $\nu_{\text{PWA}}^{\text{min}}(\rho, \dot{\rho})$, $\nu_{\text{PWA}}^{\text{max}}(\rho, \dot{\rho})$ (compare to \ref{eq:DR_MIMO_DRCd}) which cover 95\% of the feasible area.
Further details are given in the SI.

With the limits on $\nu$, the second-order dynamic ramping constraints are completely parameterized and have the form:
\begin{gather}
    \rho^{(2)} = \nu 
    \\
    a_0^{\text{min}} + a_1^{\text{min}}\rho \leq \dot{\rho} \leq a_0^{\text{max}} + a_1^{\text{max}}\rho
    \\            
    \label{eq:DR_MIMO_bounds_nu}
    \nu_{\text{PWA}}^{\text{min}} \left(\rho,\dot{\rho}\right) \leq \nu \leq \nu_{\text{PWA}}^{\text{max}}\left(\rho,\dot{\rho} \right)
\end{gather}

\subsection{Investigation 1: Ramp optimizations}

\begin{figure*}[h]
\centering 
 \includegraphics{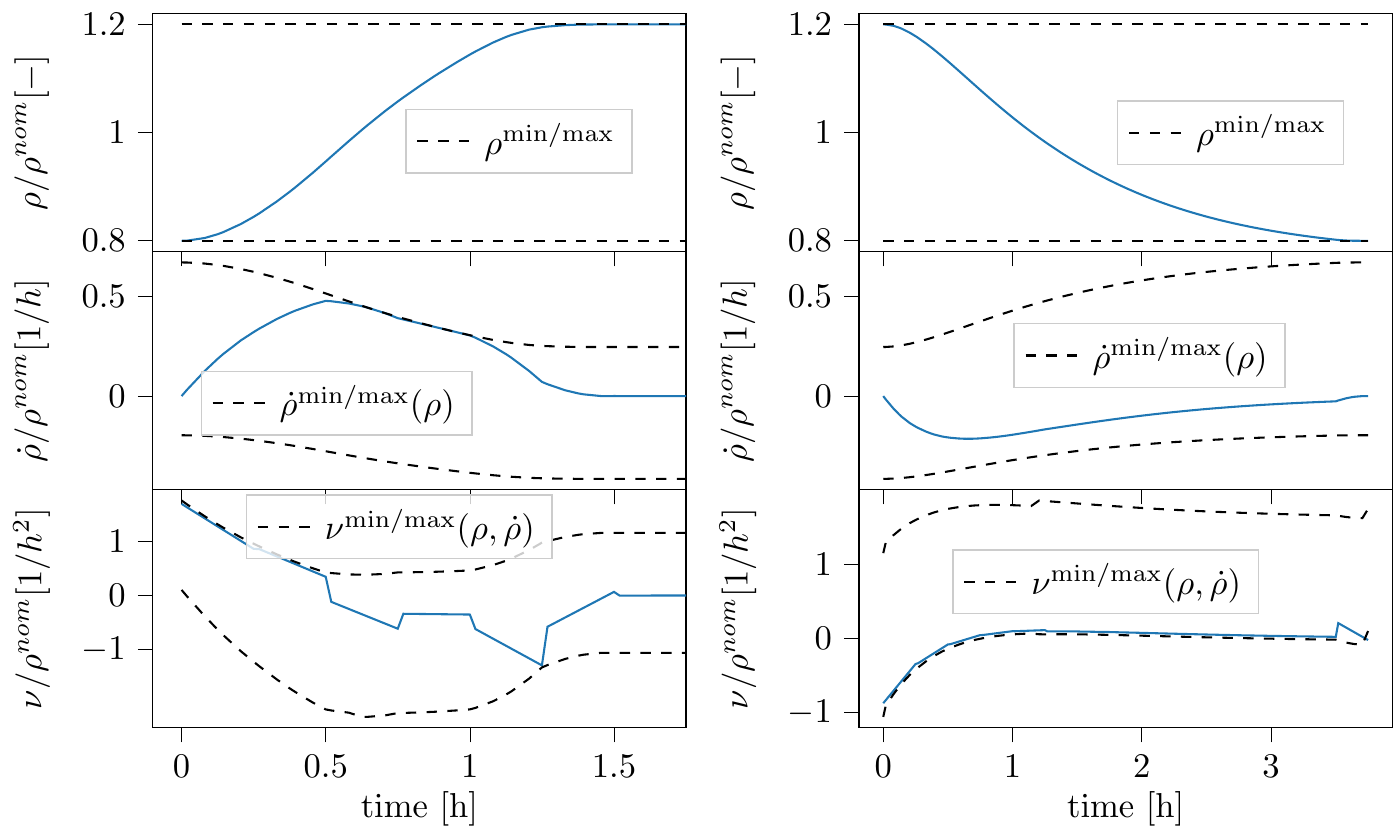}
  \caption{Production rate $\rho$, its first time derivative $\dot{\rho}$, and ramping degree of freedom $\nu$ for an as-fast-as-possible ramp up (left) and an as-fast-as-possible ramp down (right) together with their bounds. The bounds of $\dot{\rho}$ are functions of $\rho$ and the bounds of $\nu$ are functions of $\rho, \dot{\rho}$.  The ramping degree of freedom $\nu$ is discretized to be piece-wise linear.}
   \label{fig:DRMIMO_ramps}
\end{figure*}

To illustrate the ramping behavior with the derived dynamic ramping constraints, we perform two as-fast-as-possible ramp optimizations shown in Figure \ref{fig:DRMIMO_ramps}.
The ramp-up from minimum production rate to maximum production rate takes 1.3 hours, and the corresponding ramp-down takes 3.1 hours.
The ramp-up is first constrained by the acceleration, i.e., the bounds on $\nu$, and then by the speed, i.e., the bounds on $\dot{\rho}$.
In contrary, the ramp-down is always constrained by the bounds on $\nu$.

\begin{figure*}[h]
\centering
 \includegraphics[scale=0.95]{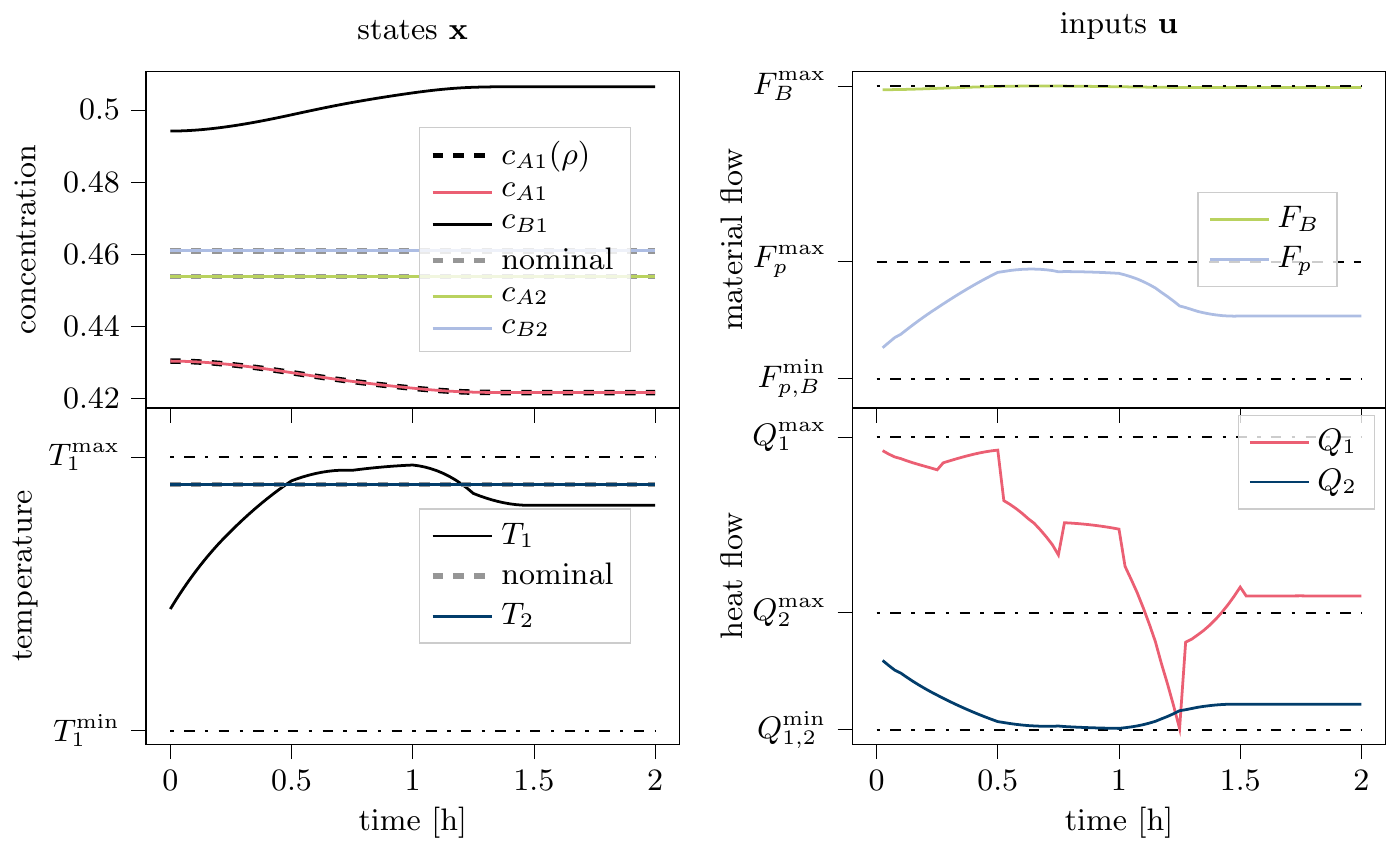}
  \caption{Simulation result of ramp-up optimization (Figure 10, left) on the full-order process model showing states $\vect{x} = (c_{A1}, c_{B1}, T_1,  c_{A2}, c_{B2}, T_2)^T$ (left), and control inputs $\vect{u} = (F_B, F_p, Q_1, Q_2)^T$ (right). Minimum (min) and maximum (max) values are shown in dashed-dotted lines. The nominal values of the flat output components $c_{A2}, c_{B2}, T_2$ and the operating strategy $\pi_4({\rho})$ for the flat output component $c_{A1}$ are shown in dashed lines.}
   \label{fig:DRMIMO_sim_ramp}
\end{figure*}

To visualize such an optimized ramp on the full-order process model, we simulate the ramp-up by using the optimized production rate trajectory (left part of Figure~\ref{fig:DRMIMO_ramps}) as input to a simulation and calculate the control inputs $\vect{u}$ using the backtransformation function $\vect{u} = \boldsymbol{\psi}_2(\rho, \dot{\rho},\nu)$ derived above.
While the first three flat output components $c_{A2}, c_{B2}, T_2$ are maintained at their nominal values, the fourth output component $c_{A1}$ follows the function of the production rate $c_{A1}=\pi_4(\rho)$ specified in the operating strategy (Figure \ref{fig:DRMIMO_sim_ramp}).
All other states and inputs are within their respective bounds. 
Moreover, in the first half hour, when the ramp-up is limited by the limit on the acceleration $\nu^{\text{max}}$ (compare to Figure \ref{fig:DRMIMO_ramps}, left), $Q_1$ is close to its maximum value (Figure \ref{fig:DRMIMO_sim_ramp}) because the upper limit of $\nu^{\text{max}}$ is derived from $Q_1^{\text{max}}$.
However, $Q_1$ does not reach its maximum value due to the conservative piecewise affine approximation.
During the second half hour, the ramp-up is limited by the maximum speed $\dot{\rho}^{\text{max}}$ (compare to Figure \ref{fig:DRMIMO_ramps}, left) and thus the control input $Q_2$, which limits the speed for high production rates $\rho$ (compare to Figure \ref{fig:DRMIMO_bounds_rho_1}), is close to its bounds (Figure \ref{fig:DRMIMO_sim_ramp}). Here, $Q_2$ comes very close to its bound as the conservative approximation of the ramping limit on $\dot{\rho}$ is very close to the true nonlinear limit (Figure \ref{fig:DRMIMO_bounds_rho_1}).
Finally, at hour 1.25, the acceleration $\nu$ touches the lower limit $\nu^{\text{min}}$ (Figure~\ref{fig:DRMIMO_ramps}) and the input $Q_1$ reaches its lower limit $Q_1^{\text{min}}$ (Figure~\ref{fig:DRMIMO_sim_ramp}).

Ramp optimization and corresponding simulation show that even though the dynamic ramping constraints are formed by linear and piecewise affine equations, they can capture dynamics that are significantly more complex than traditional static first-order ramps.
Moreover, the original process model with 6 states and 5 degrees of freedom is reduced to a dynamic ramping constraint with only 2 states, i.e., $\rho$, $\dot{\rho}$, and one degree of freedom $\nu$.
Accordingly, coupling the flat outputs to the production rate reduces the model size and thus simplifies optimization.

\subsection{Investigation 2: Demand response application}
\begin{figure*}[h]
\centering
 \includegraphics{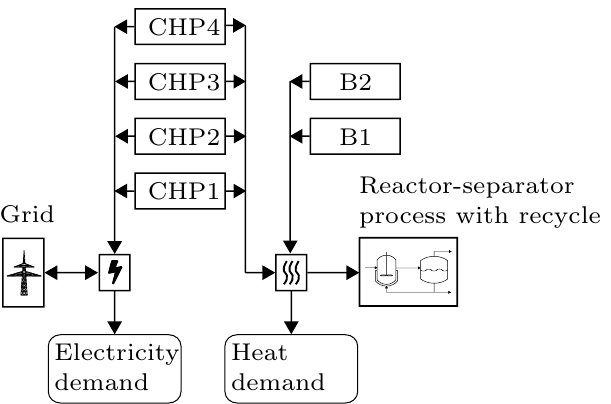}
  \caption{Multi-energy system and reactor-separator process with recycle: The flexible process as well as additional non-flexible heat and electricity demands are supplied by an multi-energy system, consisting of 4 combined heat and power plants (CHP1-CHP4) and 2 boilers (B1, B2). Moreover, electricity can be exchanged with the grid.}
   \label{fig:DRMIMO_case_Study_energy_system}
\end{figure*}
To demonstrate the dynamic ramping constraints in a DR application, the flexible process is considered together with a multi-energy system and additional non-flexible heat and electricity demands (Figure \ref{fig:DRMIMO_case_Study_energy_system}).
We consider an energy system based on \cite{SusanneSass.2020}.
However, instead of one combined heat and power plant (CHP) and one boiler (B), we extended the system to 4 CHPs and 2 boilers to study a larger energy system with more discrete on/off-decisions. 
Additionally, electricity can be bought from and sold to the grid for the day-ahead price that may change hourly.
The details about the multi-energy system and the derivation of the optimization problem (P) are given in the SI. 

The optimization problem (P) is formulated using pyomo \citep{hart2017pyomo,hart2011pyomo} and discretized using the extension pyomo.dae \citep{Nicholson.2018}.
We apply discretization by orthogonal collocation on finite elements with 2 elements per hour and 3 collocation points per element.
This discretization was found to be sufficiently accurate in preliminary calculations.
Overall, we have 144 discretization points over the 24 hour time horizon.
The discretized problem has 5,559 continuous variables. 
For the binary variables, we use a time discretization of one hour to match the time step of electricity prices.
There are 6 binaries per hour for the on/off-decision of energy system units and 3 binaries per hour for the piecewise affine limits on the ramping degree of freedom $\nu$ and the piecewise affine energy demand model (more details in the SI).
Thus, the final MILP problem has $9\times24=216$ binary variables.

The optimization problem is solved using gurobi version 9.5.2 \citep{gurobi}.
As \cite{IiroHarjunkoski.2014} state that the maximum acceptable optimization runtime for scheduling problems is typically between 5 and 20 minutes, we set the maximum optimization runtime to 5 minutes.
All calculations are performed on a Windows 11 machine with an Intel(R) Core(TM) i7-1165G7 CPU and 16 GB RAM.

\begin{figure*}[h!]
\centering
 \includegraphics[scale=0.5]{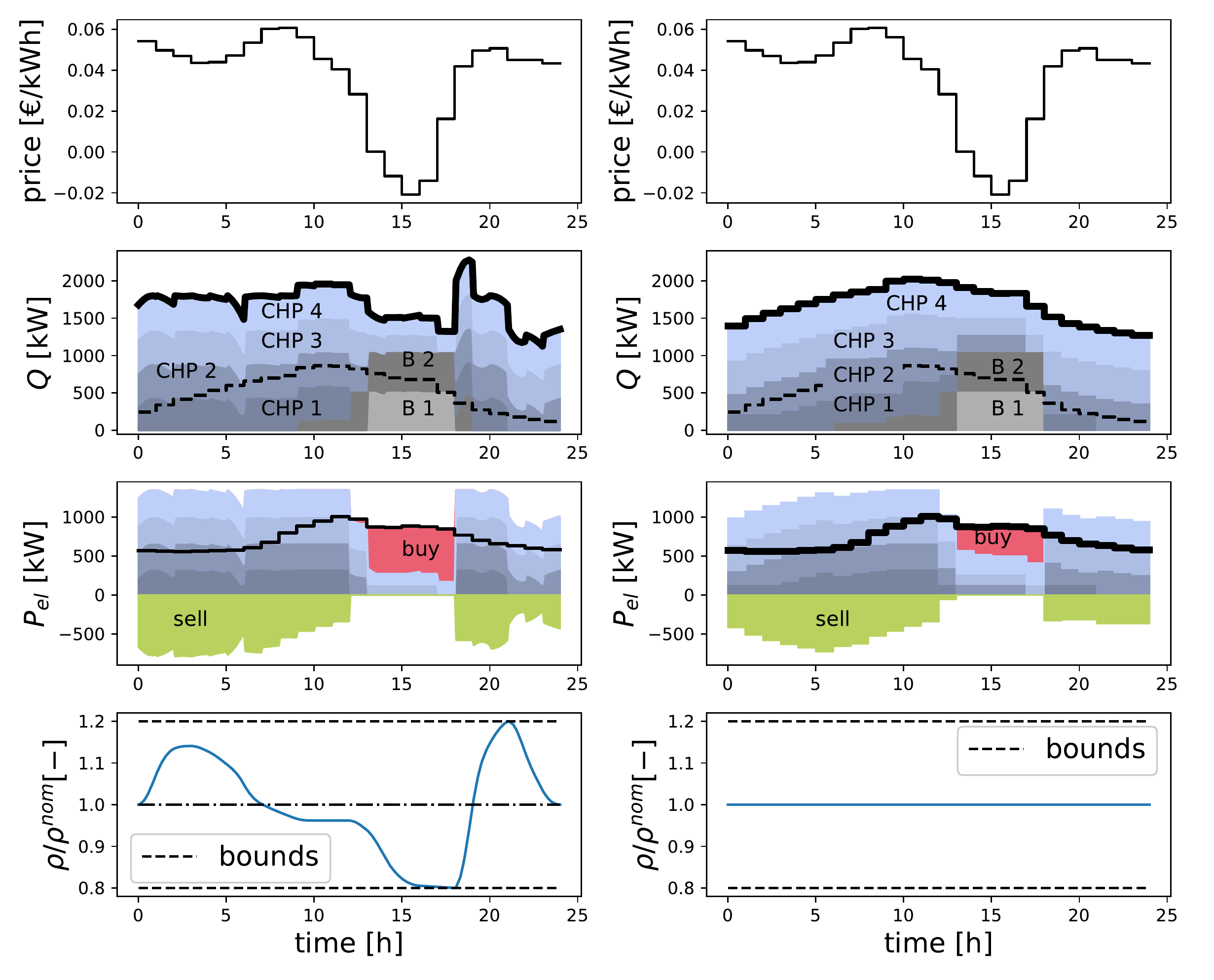}
  \caption{Resulting schedule for the flexible process performing demand response (left) and the flexible process being operated in steady-state (right). Top: electricity price. Second row: heat supplied by the 4 combined heat and power plants (CHP1 - CHP4) and the two boilers (B1, B2). The inflexible heat demand is shown as dashed black line and the total heat demand as bold black line. Third row: Electricity generated by CHP1 - CHP4, sold electricity to the grid, and electricity bought from the grid. The inflexible electricity demand is shown as black line. Bottom: Production rate $\rho$.}
   \label{fig:DRMIMO_op_DR}
\end{figure*}

First, the energy system operation is optimized without accounting for the energy demand of the flexible process to obtain the energy costs of the inflexible demands only.
Second, the operation of the energy system is optimized with the flexible process operating in steady-state such that the demand of the flexible process is constant (compare to Figure \ref{fig:DRMIMO_op_DR}, right).
Third, the DR optimization is performed using dynamic ramping constraints (compare to Figure \ref{fig:DRMIMO_op_DR}, left).

The DR operation of the process reduces the total energy costs by 4.6~\% compared to steady-state operation.
Considering only the energy costs associated with the flexible process, the cost reduction through demand response is 9.8~\%.

In the resulting operation, two types of periods can be distinguished: In times of low electricity prices, heat is preferably produced by the boilers, and electricity is bought from the grid (hours 13 - 18).
In times of high electricity prices, heat is preferable provided by the CHPs, and excess electricity is sold to the grid (hours 1-12, 19-24).
The DR case reduces costs due to two reasons:
First, the boilers are operated less.
Instead of 12 hours, the boilers are only active in 10 hours (Figure \ref{fig:DRMIMO_op_DR}).
The amount of heat provided by the boilers is reduced by 2~\%.
Second, the heat demand of the flexible process and, therefore, the electricity production of the CHPs is shifted from hours of low prices to hours with higher prices (Figures \ref{fig:DRMIMO_op_DR}).
For instance, the heat demand is lower in hour 15 and higher in hour 5.
Consequently, the derived dynamic ramping constraints allow to reduce costs substantially compared to steady-state process operation.

The optimization problem terminates after a maximum runtime of 5 minutes that we set following \cite{IiroHarjunkoski.2014}. 
The remaining optimality gap is 3.9~\%. 
As a further comparison, we simplify the bounds of the ramping degree of freedom $\nu$ in Equation~(\ref{eq:DR_MIMO_bounds_nu}) from piecewise affine to purely linear functions.
These linear functions only cover 80\% instead of 95\% of the feasible area for $\nu$ (compare to Section~\ref{subsec:DRMIMO_case_study_stratey}). 
With these linear limits $\nu_{\text{lin}}^{\text{min}}$, $\nu_{\text{lin}}^{\text{max}}$, we receive very similar results compared to those previously obtained with the piecewise affine limits: While the cost reductions achieved with the linear limits slightly  improve to 4.7\% instead of 4.6\%, the remaining optimality gap slightly worsens to 4.2\% instead of 3.9\%. 
That is, with purely linear dynamic ramping constraints, we find a comparable near-optimal solution with a comparable optimality gap.
Overall, even if the on/off-status of 6 energy system units has to be optimized simultaneously with the process operation, the optimization provides a schedule achieving substantial cost reductions within this maximum optimization runtime of 5 minutes.

\section{Discussion}
\label{sec:DRMIMO_discussion}
In this section, we discuss possible adaptations and limitations of our approach.

Throughout this paper, we assume that the production rate $\rho$ is a degree of freedom (compare to assumptions in Section~\ref{subsec:DRMIMO_assumptions}).
However, the method can be adapted to cases where the production rate is not a degree of freedom but a component of the flat output vector $\xi_k$.
For instance, if the flow rate of the product stream would be hydraulically driven by a filling level $h_{\text{fill}}$, the production rate would be given as part of the flat output by a function of $h_{\text{fill}}$.
In such cases, the production rate cannot be controlled directly but only through manipulating the corresponding input $u_k$.
For instance, the filling level $h_{\text{fill}}$ could be controlled by manipulating the feedflow of the corresponding unit. 

In our operating strategy discussed in Section~\ref{sec:DRMIMO_op_strategy}, all flat output components without given specifications are coupled to the production rate.
This coupling reduces the dimensionality of the model.
Instead of $m$ flat outputs and their derivatives, only the production rate and its derivatives are variables in the optimization problem.
Consequently, there are $\delta$ states and one degree of freedom.
This dimensionality reduction strongly reduces the computational complexity of the optimization problem.
Still, coupling all flat outputs with the production rate might be unfavorable in cases where some flat outputs can only be changed much slower than the production rate.
In such cases, the outputs without specifications could be kept as independent integrator chains in the ramping model (right part in Figure~\ref{fig:DRMIMO_steps}).
In our case study, the concentration $c_{A1}$, which is the fourth flat output component, could be uncoupled such that there are two integrator chains in the ramping model: one for the production rate~$\rho$ and one for the concentration~$c_{A1}$.
Consequently, dynamic ramping constraints need to be derived for both integrator chains.
While this uncoupled version makes the ramping model computationally more challenging, it is also more flexible and thus might enable higher profits in some cases.
An extreme example is electrolyzers that can often adapt their production rate rapidly but have slow temperature dynamics \citep{Simkoff.2020,BenjaminFlamm.2021}.
Thus, for electrolyzers, it is not favorable to couple the temperature with the production rate.
Instead, we have shown recently in a conference paper that it is favorable to keep both production rate and temperature as degrees of freedom and formulate dynamic ramping constraints on the temperature \citep{Baader.2021b}.

Conceptually, the electrolyzer example \citep{Baader.2021b} also shows that it is possible to consider two scheduling-relevant variables at a time. 
Especially, our approach based on dynamic ramping constraints could be applied to the case of processes with two production rates $\rho_1$ and $\rho_2$, where the process model from Equation~(\ref{eq:DR_MIMO_process_model}) changes to $\dot{\vect{x}} = \vect{f} (\vect{x},\vect{u},\rho_1, \rho_2)$.
In that case, there would be two ramping state vectors $\boldsymbol{\varphi}_1$ and $\boldsymbol{\varphi}_2$ and two ramping degrees of freedom $\nu_1$ and $\nu_2$. 
In general, the ramping degrees of freedom $\nu_1$, $\nu_2$ are limited by both ramping state vectors $\boldsymbol{\varphi}_1$, $\boldsymbol{\varphi}_2$ and also influence each other.
Thus, the ramping limits read:
\begin{align}
    \nu_1^{\text{min}}(\boldsymbol{\varphi}_1,\boldsymbol{\varphi}_2,\nu_2) \leq \nu_1 \leq \nu_1^{\text{max}}(\boldsymbol{\varphi}_1,\boldsymbol{\varphi}_2,\nu_2) \\
    \nu_2^{\text{min}}(\boldsymbol{\varphi}_1,\boldsymbol{\varphi}_2,\nu_1) \leq \nu_2 \leq \nu_2^{\text{max}}(\boldsymbol{\varphi}_1,\boldsymbol{\varphi}_2,\nu_1) 
\end{align}
In principle, our dynamic ramping approach can also be applied to more than two scheduling-relevant variables. However, increasing the number of scheduling-relevant variables also increases the number of arguments entering the functions $\nu_i^{\text{min}}(\cdot)$, $\nu_i^{\text{max}}(\cdot)$.
Thus, a piecewise affine approximation of the ramping limits $\nu_i^{\text{min}}(\cdot)$, $\nu_i^{\text{max}}(\cdot)$ might require many binaries and lead to a high computational burden.

The present work focuses on demand response applications where the production rate $\rho$ is the main scheduling-relevant variable.  
However, it is straightforward to adapt the approach to cases where a different variable is scheduling-relevant. 
For example, in multi-product processes, the concentration may be varied to yield different products. Thus, the scheduling needs to account for the dynamics of the concentration during product transitions \citep{FloresTlacuahuac.2006,baader2022simultaneous}. 
Our approach can be transferred to such a multi-product process if the production rate $\rho$ is replaced by the concentration.

The main limitation of our approach is the assumption of a flat process model.
For non-flat process models, the ramping limits could still be derived as functions of all process states $\vect{x}$.
However, it is not possible to find the coordinate transformation from the process states $\vect{x}$ to a transformed state vector $\boldsymbol{\Xi}$ as defined in Equation~(\ref{eq:DRMIMO_Xi}), and thus,
the ramping limits cannot be given as functions of the ramping state vector $\boldsymbol{\varphi}$ only. 
Hence, an extension to non-flat processes would be partly heuristic as the state vector $\vect{x}$ in the ramping limits would have to be approximated based on the ramping state~$\boldsymbol{\varphi}$.

Another limitation of this work is that solving the equation system to derive the backtransformation requires a lot of manual effort, even though a computer algebra system is used as support. 
It is an open question to what extent our approach can be automated to enable the analysis of large-scale processes. 
Still, we do not consider this to be a significant restriction of our work as processes fulfilling our flatness assumption usually do not have too many states \citep{OLDENBURG2002385}:
As inputs must be matched to outputs such that all states are covered (compare to the necessary graphical condition, e.g., Figure~\ref{fig:DRMIMO_graph_case_study}), processes having much more states than inputs are usually not flat. 
Thus, our method will likely not be applicable to large-scale process models with hundreds of states. 
However, our approach seems applicable for flat process models with a number of states in the low double-digit range.
As stated at the beginning of Section~\ref{sec:DRMIMO_case_study}, we only reduced the textbook example by \cite{Christofides.2011} from 9 states to 6 to ease the readability of the paper.
Still, in a previous version, we considered the original version with 9 states and could find the inverse transformation easily using computer algebra.
In future work, our approach could be coupled with model-order reduction approaches that derive a low-order representation of the process dynamics for the slow time scale \citep{baldea_daoutidis_2012}. 
Even if the full-scale process model is not flat, the low-order dynamics relevant for scheduling optimization might be.

\section{Conclusion}
\label{sec:DRMIMO_conlusion}
Dynamic ramping constraints simplify the simultaneous demand response scheduling optimization of production processes and their energy systems compared to an optimization considering the full-order process model. 
Still, dynamic ramping constraints can capture more flexibility than traditional static ramping constraints as they allow high-order dynamics and non-constant ramp limits.
In this paper, we extend our method to rigorously derive dynamic ramping constraints  from input-state linearizable single-input single-output (SISO) processes \citep{Baader.2021} to flat multi-input multi-output (MIMO) processes. 
In the MIMO case, dynamic ramping constraints reduce the problem dimensionality by coupling all flat outputs with the production rate.
In our case study, a system with 6 states and 5 degrees of freedom is reduced to 2 states and one degree of freedom.
Additionally, we demonstrate that an operational strategy can be chosen for flat outputs such that the steady-state production points are optimal, e.g., with respect to energy consumption.

Our case study demonstrates that dynamic ramping constraints allow for finding DR schedules for flexible processes and multi-energy systems that substantially reduce energy costs compared to a steady-state operation.
Even though discrete on/off-decisions in the multi-energy system add to the computational complexity, the problem can be solved within the time limit for online scheduling.

Overall, dynamic ramping constraints allow bridging the gap between nonlinear process models and simplified process representations suitable for real-time scheduling optimization.
  
  \section*{Author contributions}
\textbf{Florian J. Baader (FB)}: Conceptualization, Methodology, Software, Investigation, Validation, Visualization, Writing - original draft. 
\textbf{Philipp Althaus (PA)}: Conceptualization, Writing – review \& editing.
\textbf{André Bardow (AB)}: Funding acquisition, Conceptualization, Supervision, Writing – review \& editing. 
\textbf{Manuel Dahmen (MD)}: Conceptualization, Supervision, Writing – review \& editing.

  \section*{Declaration of Competing Interest}
    We have no conflict of interest.

    \section*{Acknowledgements}
    This work was supported by the Helmholtz Association, Germany under the Joint Initiative ‘Energy System Integration’. AB and FB also received support from the Swiss Federal Office of Energy, Switzerland through the project ‘SWEET PATHFNDR’.

\newpage
\twocolumn
\section*{Nomenclature} \label{sec:nomenclature}

\noindent\textbf{Abbreviations} \\
\noindent
\begin{tabularx}{\columnwidth}{lX}
CHP & combined heat and power plant \\
CSTR & continuous stirred tank reactor\\
DR & demand response \\
DRC & dynamic ramping constraint \\
MIDO & mixed-integer dynamic optimization \\
MILP & mixed-integer linear program \\
MIMO & multi-input multi-output \\
MINLP & mixed-integer nonlinear program \\
SISO & single-input single-output 
\end{tabularx} \\

\noindent\textbf{Greek symbols} \\
\noindent
\begin{tabularx}{\columnwidth}{lX}
$\alpha$ & integer number (compare to assumption 3a)\\
$\alpha_K$ & relative volatility of component K \\
$\beta$ & integer number  (compare to assumption 3b)\\
$\gamma$ & integer number (compare to assumption 3b)\\
$\delta$ & order of dynamic ramping constraint \\
$\zeta$ & integer number \\
$\kappa$ &  integer number (compare to assumption 3a)\\
$\nu$ & ramping degree of freedom \\
$\boldsymbol{\Xi}$ & transformed state vector \\
$\xi$ & flat output \\
$\pi$ & operating strategy \\
$\rho$ & production rate \\
$\varrho_F$ & density \\
$\Phi$ & objective \\
$\phi$ & transformation to flat output \\
$\boldsymbol{\varphi}$ & ramping state vector \\
$\psi$ & backtransformation from flat output \\
\end{tabularx} \\

\noindent\textbf{Latin symbols} \\
\noindent
\begin{tabularx}{\columnwidth}{lX}
$a$ & coefficient \\
$b$ & parameter \\
$C$ & coverage \\
$C_p$ & heat capacity \\
$c$ &concentration \\
$E$ & activation energy \\
F & flow rate \\
$f$ & nonlinear function \\
$H$ & enthalpy \\
$k$ & reaction constant \\
$m$ & number of control inputs \\
$n$ & number of states \\
$P$ & power \\
$Q$ & heat flow \\
$R$ & gas constant \\
$T$ & temperature \\
$t$ & time \\
$u$ & control input \\
\end{tabularx} 
\begin{tabularx}{\columnwidth}{lX}
$V$ & volume \\
$v$ & vertex \\
$x$ & state \\
$z$ & binary variable \\
\end{tabularx}

\noindent\textbf{Subscripts} \\
\noindent
\begin{tabularx}{\columnwidth}{lX}
0 & feed \\
1 & reactor \\
2 & flash \\
A & component A \\
B & component B or bottom \\
C & component C \\
dem & demand \\
el & electric \\
f & final time  \\
p & purge \\
s & steady-state \\
V & vaporization \\
\end{tabularx} \\

\noindent\textbf{Superscripts} \\
\noindent
\begin{tabularx}{\columnwidth}{lX}
max & maximum \\
min & minimum \\
nom & nominal \\
\end{tabularx} \\
\onecolumn

\bibliographystyle{apalike}
\renewcommand{\refname}{Bibliography} 
\bibliography{literature.bib}

\end{document}


\thispagestyle{firststyle}

  \begin{center}
    \begin{large}
    {\fontsize{12}{14} \selectfont
      \textbf{\mytitle}}
    \end{large} \\
    \myauthor
  \end{center}

  \begin{footnotesize}
    \affil
  \end{footnotesize}

\section{Exemplary model for necessary flatness conditions}
\begin{figure*}[h]
\centering
\includegraphics{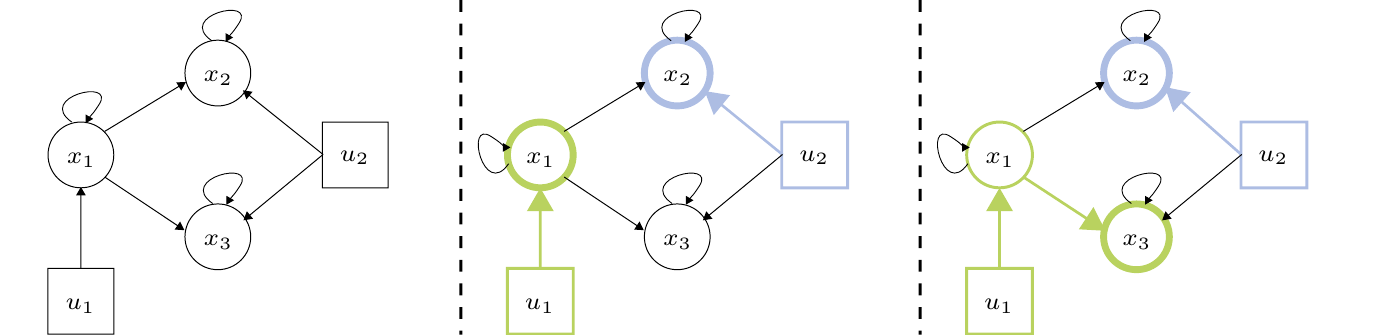}
  \caption{Graph representation of example process (E) (left). While the candidate  $\boldsymbol{\xi}=\left(x_1, x_2 \right)^T$ does not fulfil the necessary condition for a flat output vector (center), the candidate  $\boldsymbol{\xi}=\left(x_3, x_2\right)^T$ does (right).}
   \label{fig:DRMIMOSI_graph_representation}
\end{figure*}

In this section, the two necessary conditions for flatness (compare to Section 3.2.1 in the main paper) are illustrated using an example (E) with three states and two inputs:
    \begin{align}
        &\dot{x}_1 = -x_1 - u_1 \tag{Ea}\\
        &\dot{x}_2 = -x_2 + x_1 + u_2 \tag{Eb}\\
        &\dot{x}_3 = -x_3 + 2 x_1 - u_2 \tag{Ec}
    \end{align}
In the example (E), the output vector $\boldsymbol{\xi}=\left( x_1, x_2  \right)^T$ does not fulfil the necessary graphical condition for flat outputs because the state $x_3$ is not covered (Figure \ref{fig:DRMIMOSI_graph_representation}, center).
However, the vector $\boldsymbol{\xi}=\left( x_3,  x_2 \right)^T$ is a candidate for a flat output as all states are covered (Figure \ref{fig:DRMIMOSI_graph_representation}, right).

To receive the equation system for the backtransformations, the flat output candidate $\boldsymbol{\xi}=\boldsymbol{\phi} \left( \vect{x} \right)=\left(\begin{array}{c} x_3 \\ x_2 \end{array} \right)$ has to be differentiated two times to obtain a square backtransformation equation system:
\begin{align}
    \boldsymbol{\xi}&=\boldsymbol{\phi} \left( \vect{x} \right)=\left(\begin{array}{c} x_3 \\ x_2 \end{array} \right) \label{eq:DRMIMOSI_backtr_1}\\
    \dot{\boldsymbol{\xi}} &= \frac{\partial \boldsymbol{\phi}(\vect{x})}{\partial \vect{x}} \dot{\vect{x}}  = \frac{\partial\left(\begin{array}{c} x_3 \\ x_2 \end{array} \right) }{\partial \vect{x}} \vect{f}(\vect{x},\vect{u})\\ 
    \boldsymbol{\xi}^{(2)} &= \frac{\partial \dot{\boldsymbol{\xi}}}{\partial \vect{x}} \dot{\vect{x}} + \frac{\partial \dot{\boldsymbol{\xi}}}{\partial \vect{u}} \dot{\vect{u}} \label{eq:DRMIMOSI_backtr_3}
\end{align}
Structurally, the backtransformation equation system (\ref{eq:DRMIMOSI_backtr_1}) - (\ref{eq:DRMIMOSI_backtr_3}) is solvable for the three states $x_1, x_2, x_3$, the two inputs $u_1, u_2$, and the input derivative $\dot{u}_2$ (Table \ref{tab:DRMIMOSI_incidence_E}).
The input derivative $\dot{u}_2$ is part of the equation system as both second derivatives $\xi_1^{(2)}$, $\xi_2^{(2)}$ are functions of $\dot{u}_2$.

\begin{table}[bt]
    \centering
    \setlength{\tabcolsep}{3pt}
    \renewcommand{\arraystretch}{1.0}
    \caption{Analysis of structural solvability of the backtransformation equation system to calculate states $\vect{x}$, and inputs $\vect{u}$ with flat output candidate $\boldsymbol{\xi}=\left( x_3,  x_2  \right)^T$ for example process (E). All states and inputs that appear in an equation are marked with a cross. A necessary condition for solvability is that it must be possible to circle exactly one cross in every row and exactly one cross in every line. In other words: Each equation should be solved for one variable.}
    \label{tab:DRMIMOSI_incidence_E}
    \begin{tabular}{|c |p{0.7cm}|p{0.7cm}|p{0.7cm}|p{0.7cm}|p{0.7cm}|p{0.7cm}|}
    &\multicolumn{2}{c}{states} &&\multicolumn{3}{p{2.1cm}|}{inputs and derivatives}\\
    &$x_{1}$& $x_{2}$& $x_{3}$& $u_{1}$& $u_{2}$ &$\dot{u}_2$\\
    \hline
            $\xi_1$ & &  &\circled{x}  &  &  & \\
            $\xi_2$ & &\circled{x}  &  &  &  & \\
            $\dot{\xi}_1$ &\circled{x} &  &x  &  &x  & \\
            $\dot{\xi}_2$ &x &x  &  &  &\circled{x}  & \\
            $\xi_1^{(2)}$ &x &  &x  &\circled{x}  &x  &x \\
            $\xi_2^{(2)}$ &x &x  &  &x  &x  & \circled{x}\\
    \end{tabular}
\end{table}

\section{Illustrative example for multiple feasible regions}
As discussed in Section 3.4 of the main manuscript, multiple distinct feasible regions can occur for the dynamic ramping constraints.
To illustrate this phenomenon, we consider an illustrative model (IM) with three states $\vect{x}$,  two inputs $\vect{u}$, and a production rate $\rho$:
    \begin{align}
        \dot{x}_1 &= x_2 + \rho \tag{IMa} \label{eq:DRMIMOSI_IPa}\\
        \dot{x}_2 &= -x_2 + u_1 \tag{IMb} \label{eq:DRMIMOSI_IPb}\\
        \dot{x}_3 &= x_2^2 + u_2 \tag{IMc}  \label{eq:DRMIMOSI_IPc}
    \end{align}
As shown in the following, this process has a flat output $\boldsymbol{\xi} = \boldsymbol{\phi} \left(\vect{x}\right) =\left( \begin{array}{c} x_1 \\ x_3 \end{array} \right)$.
We assume an operating strategy $\boldsymbol{\xi} = \boldsymbol{\pi}(\rho)=\left( \begin{array}{c} a \rho \\ b \rho \end{array} \right)$ with two constants $a$, $b$.

For $x_1$, and $x_3$, the back transformation is given by the function $\boldsymbol{\pi}$.
 For $x_2$, Equation~(\ref{eq:DRMIMOSI_IPa}) is set equal to $\dot{\xi}_1 = \dot{x}_1 = a\dot{\rho}$ which gives
    \begin{align}
        x_2 = \underbrace{a\dot{\rho}}_{\dot{x}_1} - \rho
    \end{align}
Similarly, inserting $\dot{x}_2 = a\rho^{(2)} - \dot{\rho}$ into Equation~(\ref{eq:DRMIMOSI_IPb}) gives
    \begin{align}
        u_1 = \underbrace{a\dot{\rho} - \rho}_{x_2} + \underbrace{a\rho^{(2)} - \dot{\rho}}_{\dot{x}_2 }
    \end{align}
Finally, inserting $\dot{\xi}_2 = \dot{x}_3 = b\dot{\rho}$ into Equation~(\ref{eq:DRMIMOSI_IPc}) gives
    \begin{align}
        \label{eq:DRMIMOSI_u_2}
        u_2 = \underbrace{b\dot{\rho}}_{\dot{x}_3} - \left( \underbrace{a\dot{\rho} - \rho}_{{x_2}} \right) ^2
    \end{align}
As all states and inputs can be given as functions of production rate $\rho$ and derivatives $\dot{\rho}$, $\rho^{(2)}$, the process is de facto flat for the given operating strategy.

\begin{figure}[h]
\centering
 \includegraphics{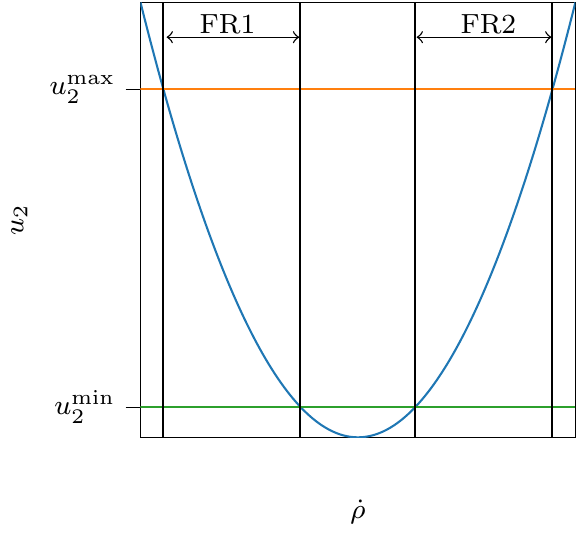}
  \caption{Sketch for the input $u_2$ as function of the first derivative of the production rate $\dot{\rho}$ (compare to Equation~(\ref{eq:DRMIMOSI_u_2})). Potentially, there can be two feasible regions for $\dot{\rho}$, FR1 and FR2, respectively.}
   \label{fig:DRMIMOSI_concept_two_regions}
\end{figure}

Two feasible regions for the dynamic ramping constraints can result from the second input $u_2$ being a quadratic function of the first derivative of the production rate $\dot{\rho}$ (Equation~(\ref{eq:DRMIMOSI_u_2})) (Figure \ref{fig:DRMIMOSI_concept_two_regions}).
One could either restrict the operating range to one of the two regions or would have to introduce a binary variable that indicates which region is active.

\section{Supplementary information to case study}

\subsection{Model equations}
In this section, we present the model equations of the case study (compare to Figure 4 in the main paper). 
All symbols other than states $\vect{x} = (c_{A1}, c_{B1}, T_1,  c_{A2}, c_{B2}, T_2)^T$, control inputs $\vect{u} = (F_B, F_p, Q_1, Q_2)^T$, and production rate $\rho$ are parameters marked in gray with values listed in Table \ref{tab:DRMIMOSI_parameters_1}.
The index $0$ indicates the feed quantities.
\begin{align}
    \label{eq:DRMIMOSI_dc_A1}
    \dot{c}_{A1} &= \frac{\rho + F_p}{\textcolor{gray}{V_1}}(\textcolor{gray}{c_{A0}} - c_{A1}) + \frac{F_B - F_p}{\textcolor{gray}{V_1}}(c_{A2} - c_{A1}) - \textcolor{gray}{k_1} c_{A1}\textcolor{gray}{e}^{-\frac{\textcolor{gray}{E_1}}{\textcolor{gray}{R}T_1}}
    \\
    \label{eq:DRMIMOSI_dc_B1}
    \dot{c}_{B1} &= \frac{\rho + F_p}{\textcolor{gray}{V_1}}(\textcolor{gray}{c_{B0}} - c_{B1}) + \frac{F_B - F_p}{\textcolor{gray}{V_1}}(c_{B2} - c_{B1}) 
    \\ \nonumber 
    &+ \textcolor{gray}{k_1} c_{A1}\textcolor{gray}{e}^{-\frac{\textcolor{gray}{E_1}}{\textcolor{gray}{R}T_1}} -  \textcolor{gray}{k_2} c_{B1}\textcolor{gray}{e}^{-\frac{\textcolor{gray}{E_2}}{\textcolor{gray}{R}T_1}}
    \\
    \label{eq:DRMIMOSI_dT_1}
   \dot{T}_{1} &= \frac{\rho + F_p}{\textcolor{gray}{V_1}}(\textcolor{gray}{T_{0}} - T_1) + \frac{F_B - F_p}{\textcolor{gray}{V_1}}(T_2 - T_1) - \frac{\textcolor{gray}{\Delta H_1}}{\textcolor{gray}{C_p}} \textcolor{gray}{k_1} c_{A1}\textcolor{gray}{e}^{-\frac{\textcolor{gray}{E_1}}{\textcolor{gray}{R}T_1}} 
    \\ \nonumber 
    &- \frac{\textcolor{gray}{\Delta H_2}}{\textcolor{gray}{C_p}}\textcolor{gray}{k_2} c_{B1}\textcolor{gray}{e}^{-\frac{\textcolor{gray}{E_2}}{\textcolor{gray}{R}T_1}} + \frac{Q_1}{\textcolor{gray}{\varrho_F C_p V_1}}
\end{align}

\begin{align}
    \label{eq:DRMIMOSI_dc_A2}
    \dot{c}_{A2} &= \frac{\rho + F_B}{\textcolor{gray}{V_2}}(c_{A1}-c_{A2}) -
    \\ \nonumber 
    &\frac{\rho}{\textcolor{gray}{V_2}}\left(\frac{\textcolor{gray}{\alpha_A} c_{A2}}{\textcolor{gray}{\alpha_A} c_{A2} + \textcolor{gray}{\alpha_B} c_{B2} + \textcolor{gray}{\alpha_C} (1-c_{A2}-c_{B2})}-c_{A2}\right)\\
    \label{eq:DRMIMOSI_dc_B2}
    \dot{c}_{B2} &= \frac{\rho + F_B}{\textcolor{gray}{V_2}}(c_{B1}-c_{B2}) \\ \nonumber 
    &- \frac{\rho}{\textcolor{gray}{V_2}}\left(\frac{\textcolor{gray}{\alpha_B} c_{B2}}{\textcolor{gray}{\alpha_A} c_{A2} + \textcolor{gray}{\alpha_B} c_{B2} + \textcolor{gray}{\alpha_C} (1-c_{A2}-c_{B2})}-c_{B2}\right)\\
    \label{eq:DRMIMOSI_dT_2}
    \dot{T}_{2} &= \frac{\rho + F_B}{\textcolor{gray}{V_2}}(T_1 - T_2) - \textcolor{gray}{\Delta H_V}\frac{\rho}{\textcolor{gray}{\varrho_F C_p V_2}} + \frac{Q_2}{\textcolor{gray}{\varrho_F C_p V_2}}
\end{align}

\begin{table}[h!]
    \centering
    \caption{Parameters of the reactor-separator process with recycle.}
    \label{tab:DRMIMOSI_parameters_1}
    \begin{tabular}{llr} \\
     &symbol & value  \\
     \hline
    volume CSTR & $V_1$ & 30 m³  \\
    feed concentration A & $c_{A0}$ & 1  \\
    reaction constant reaction 1& $k_1$ & 3.21$\times 10^6 \frac{1}{\text{h}}$  \\
    activation energy reaction 1 & $E_1$ & 5$\times 10^4 \frac{\text{kJ}}{\text{kmol}}$ \\
    gas constant & $R$ & 8.314  $\frac{\text{kJ}}{\text{kmolK}}$\\
    feed concentration B & $c_{B0}$ & 0 \\
    reaction constant reaction 2& $k_2$ & 3.21$\times 10^6 \frac{1}{\text{h}}$ \\
    activation energy reaction 2 & $E_2$ & 6$\times 10^4 \frac{\text{kJ}}{\text{kmol}}$  \\
    feed temperature & $T_0$ & 300 K \\
    reaction enthalpy 1 & $\Delta H_1$ & -261 $\frac{\text{kJ}}{\text{kg}}$ \\
    heat capacity & $C_p$ & 4.2 $\frac{\text{kJ}}{\text{kgK}}$ \\
    reaction enthalpy 2 & $\Delta H_2$ & -304 $\frac{\text{kJ}}{\text{kg}}$ \\
    density & $\varrho_F$ & 1000 $\frac{\text{kg}}{\text{m}^3}$  \\
    volume flash & $V_2$ & 30 m³  \\
    relative volatility of A& $\alpha_A$ & 0.5  \\
    relative volatility of B& $\alpha_B$ & 0.25  \\
    relative volatility of C& $\alpha_C$ & 1  \\
    enthalpy of vaporization& $\Delta H_V$ & 7.2$\times 10^4 \frac{\text{kJ}}{\text{kg}}$  \\
    \hline  
    \end{tabular}
\end{table}

\begin{table}[bt]
\centering
    \setlength{\tabcolsep}{3pt}
    \renewcommand{\arraystretch}{1.0}
    \caption{Analysis of structural solvability with flat output candidate $\boldsymbol{\xi}=\left( c_{A2},  c_{B2}, T_2, c_{A2}  \right)^T$ for the reactor-separator cycle (compare to Table \ref{tab:DRMIMOSI_incidence_E}).}
    \begin{tabular}{c|c|c|c|c|c|c|c|c|c|c|c|c|c}
        \label{tab:DRMIMOSI_case_sudy_incidence}
    &\multicolumn{5}{c}{states} &&\multicolumn{7}{c}{control inputs and derivatives}\\
    &$c_{A1}$& $c_{B1}$& $T_1$ & $c_{A2}$& $c_{B2}$& $T_2$ &$F_B$ &$F_p$ & $Q_1$& $Q_2$&$\dot{F}_B$&$\dot{F}_p$& $F_B^{(2)}$\\
    \hline
            $\xi_1$ & &  &  &\circled{x}   &  &  &  &  &  &  &  &  &    \\
            $\xi_2$ &  & &  &  &\circled{x}   &  &  &  &  &  &  &  &    \\
            $\xi_3$ &  &  &  & &  &\circled{x}   &  &  &  &  &  &  &     \\
            $\xi_4$ &\circled{x}   &  &  &  &  &  & &  &  &  &  &  &    \\
    \hline
    \hline
        $\dot{\xi}_1$ &x &  &  &x &x &  &\circled{x} &  &  &  &  &  &    \\
        $\dot{\xi}_2$ &  &\circled{x} &  &x &x &  &x &  &  &  &  &  &    \\
        $\dot{\xi}_3$ &  &  &x &  &  &x &x &  &  &\circled{x} &  &  &     \\
       $\dot{\xi}_4$  &x &  &\circled{x} &x &  &  &x &x &  &  &  &  &    \\
    \hline
    \hline
    $\xi_1^{(2)}$ &x &x &x &x &x &  &x &x &  &  &\circled{x} &  &    \\
    $\xi_2^{(2)}$ &x  &x &x &x &x &  &x &\circled{x} &  &  &x &  &    \\
    $\xi_4^{(2)}$ &x &x &x &x &x &x &x &x &\circled{x} &  &x &x &     \\
    \hline
    \hline
    $\xi_1^{(3)}$ &x &x &x &x &x & x&x &x &x &  &x &x &\circled{x}   \\
    $\xi_2^{(3)}$ &x &x &x &x &x &x &x &x &x &  &x &\circled{x} &x   \\
    \end{tabular}
\end{table}

\subsection{Operating strategy}
As discussed in the main paper, the fourth flat output $\xi_4$ can be chosen freely.
Here, we calculate the optimal steady-state operating points for $\xi_4$ and fix $\xi_4$ to be a function of the production rate $\xi_4=\pi_4(\rho)$.
Accordingly, we minimize the total heat input to the process $Q_1 + Q_2$ as the objective function to obtain the steady-state operating points. 
First, we determine the optimal steady-state operating points.
To this end, a steady-state optimization problem ($P_s$) is formulated by sampling the production range $\rho^{\text{min}} \leq \rho  \leq \rho^{\text{max}}$ with 21 equally distributed points $\rho_j$.
The optimization problem ($P_s$) reads:
\begin{align}
    \label{eq:DRMIMOSI_opt_ss_obj}
    &\hspace{0.4cm}  \underset{\vect{x}_j,\vect{u}_j}{\text{min}}~~  \sum_{j=1}^{21} Q_{1,j} +  Q_{2,j} \tag{$P_sa$}\\
    \text{s.t.  }&\boldsymbol{0} = \vect{f}(\vect{x}_j, \vect{u}_j, \rho_j) &~\forall j \label{eq:DRMIMOSI_opt_ss_1} \tag{$P_sb$}\\
    &\phi_{k,j}(\vect{x})  = \xi_k^{\text{nom}} &~\forall k = 1,2,3~ \forall j \label{eq:DRMIMOSI_opt_ss_2} \tag{$P_sc$}\\
    &\vect{u}^{\text{min}} \leq \vect{u}_j \leq \vect{u}^{\text{max}} &~\forall j \label{eq:DRMIMOSI_opt_ss_4}\tag{$P_sd$}\\
    &\vect{x}^{\text{min}} \leq \vect{x}_j \leq \vect{x}^{\text{max}} &~\forall j  \label{eq:DRMIMOSI_opt_ss_5} \tag{$P_se$}
\end{align}
For every point $j$, the process needs to be in steady state (Equation~(\ref{eq:DRMIMOSI_opt_ss_1})), the first three flat output components $\xi_k = \phi_{k}(\vect{x})$ need to be at their nominal values (Equation~(\ref{eq:DRMIMOSI_opt_ss_2})), and inputs and states needs to be within bounds (Constraints (\ref{eq:DRMIMOSI_opt_ss_4}), (\ref{eq:DRMIMOSI_opt_ss_5})).
All bounds are given in Table \ref{tab:DRMIMOSI_bounds}.
We implement the optimization problem $P_s$ using pyomo \citep{hart2017pyomo} and solve it using BARON version 20.10.16 \citep{Khajavirad.2018} in heuristic mode.

\begin{table}[h!]
    \centering
    \caption{State and input bounds}
    \label{tab:DRMIMOSI_bounds}
    \begin{tabular}{lrr}
    variable & lower bound & upper bound  \\
     \hline
    $c_{A1}$ & 0  & 1 \\
    $c_{B1}$ & 0 & 1\\
    $T_{1}$ & 410 K &  460K\\
    $F_{B}$ & 0$\frac{m^3}{h}$ &  20$\frac{m^3}{h}$\\
    $F_{p}$ & 0$\frac{m^3}{h}$ &  8$\frac{m^3}{h}$\\
    $Q_{1}$ & 0$\frac{kJ}{h}$  & 10$\times10^6\frac{kJ}{h}$\\
    $Q_{2}$ & 0$\frac{kJ}{h}$  & 4$\times10^6\frac{kJ}{h}$\\
    $\rho$ & 4.2$\frac{m^3}{h}$ &  6.3$\frac{m^3}{h}$\\
    \end{tabular}
\end{table}

Next, a function $\pi_4(\rho)$ needs to be chosen that does not increase the objective function (\ref{eq:DRMIMOSI_opt_ss_obj}) too much.
In other words, the steady-state operation points reached with the additional constraint
\begin{align}
    \xi_{4,j}  = \pi_4(\rho_j) &~~~~~\forall j \label{eq:DRMIMOSI_opt_ss_3} \tag{$P_sf$}
\end{align}
should give an objective function close to the objective function reached without the additional constraint (\ref{eq:DRMIMOSI_opt_ss_3}).
As a first candidate, we investigate the simplest possible operating strategy, that is to hold $\xi_4$ constant, i.e., $\pi_4(\rho) = \xi_4^{\text{nom}}$.
This optimization leads to a feasible solution. 
Accordingly, it is feasible to always hold the fourth flat output component at a constant value.
However, the objective function worsens by 2.9~\% compared to the first optimization without constraint (\ref{eq:DRMIMOSI_opt_ss_3}) because it is not possible to choose a single constant value $\xi_4^{\text{nom}}$ that is optimal for all operating points.
As a second candidate, we therefore investigate the linear strategy
\begin{align}
    \label{eq:DRMIMOSI_case_study_lin_op_strategy}
    \xi_4 = \pi_4(\rho) = a_0^{\xi_4} + a_1^{\xi_4}\rho ,
\end{align}
with optimization degrees of freedom $a_0^{\xi_4}$, $a_1^{\xi_4}$.
The optimization finds values $a_0^{\xi_4}$, $a_1^{\xi_4}$ such that the objective value worsens only 0.05~\% compared to the first optimization without the coupling constraint (\ref{eq:DRMIMOSI_opt_ss_3}). 
As this linear strategy is very close to the objective value found without constraint (\ref{eq:DRMIMOSI_opt_ss_3}), we do not consider other functions, e.g., higher-order polynomials, and continue with the linear operating strategy (\ref{eq:DRMIMOSI_case_study_lin_op_strategy}).

\subsection{Derivation of flat process reformulation}
In this section, the nonlinear coordinate transformation is derived based on the nonlinear system of equations and the operating strategy. 
The aim is to express the states $\vect{x} = (c_{A1}, c_{B1}, T_1,  c_{A2}, c_{B2}, T_2)^T$ and control inputs $\vect{u} = (F_B, F_p, Q_1, Q_2)^T$ as functions of production rate $\rho$, and its first two derivatives $\dot{\rho}$, and $\nu$.

For the reformulation, first, the strategy for the 4 flat output components $\boldsymbol{\xi} = \left(c_{A2}, c_{B2}, T_2, c_{A1}\right)^T$ discussed in Section~4.2 of the main manuscript is inserted:
    \begin{align}
        &c_{A2} = c_{A2}^{\text{nom}}~
        \\
        &c_{B2} = c_{B2}^{\text{nom}}~
        \\
        &T_{2} = T_{2}^{\text{nom}}
        \\
        \label{eq:DRMIMOSI_c_A2_rho}
        &c_{A1}(\rho) = \textcolor{gray}{a_0^{\xi_4}} + \textcolor{gray}{a_1^{\xi_4}}\rho~~ (\text{Equation (\ref{eq:DRMIMOSI_case_study_lin_op_strategy}}))
    \end{align}
While the first three outputs are chosen independent of the production rate, the fourth output $c_{A1}$ is a function of the production rate $\rho$.

Next, we consider the equation for the derivative of the first flat output, i.e., $\dot{\xi}_1 = \dot{c}_{A2}$.
As $c_{A2}$ is held constant, the derivative $\dot{c}_{A2}$ is zero and we solve Equation~(\ref{eq:DRMIMOSI_dc_A2})~for the bottom flow $F_B$ that is a function of the production rate $\rho$, i.e.,
    \begin{align}
        \label{eq:DRMIMOSI_F_B_of_rho}
        &F_B(\rho) = \rho\left( \frac{ c_{Av}^{\text{nom}} - c_{A2}^{\text{nom}}}{c_{A1}(\rho) - c_{A2}^{\text{nom}}} -1 \right),
        \\ \nonumber
        &\text{with } c_{Av}^{\text{nom}} = \frac{\textcolor{gray}{\alpha_{A}} c_{A2}^{\text{nom}}}{\textcolor{gray}{\alpha_{A}} c_{A2}^{\text{nom}} + \textcolor{gray}{\alpha_{B}} c_{B2}^{\text{nom}} + \textcolor{gray}{\alpha_{C}} \left(1- c_{A2}^{\text{nom}} - c_{B2}^{\text{nom}} \right)},
    \end{align}
    and $c_{A1}(\rho)$ as defined in Equation~(\ref{eq:DRMIMOSI_c_A2_rho}).

The second output is maintained constant, too.
Accordingly, we receive the condition $\dot{\xi}_2 = \dot{c}_{B2} = 0$.
After inserting $F_B(\rho)$ from Equation~(\ref{eq:DRMIMOSI_F_B_of_rho}) into Equation~(\ref{eq:DRMIMOSI_dc_B2}), we get the equation:
\begin{align}
    0=\frac{\rho \left(\left(c_{A1}(\rho) - c_{A2}^{\text{nom}}\right) \left(c_{B2}^{\text{nom}} - c_{Bv}^{\text{nom}}\right) + \left(- c_{A2}^{\text{nom}} + c_{Av}^{\text{nom}}\right) \left(c_{B1} - c_{B2}^{\text{nom}}\right)\right)}{\textcolor{gray}{V_{2}} \left(c_{A1}(\rho) - c_{A2}^{\text{nom}}\right)}
    \\ \nonumber
    \text{with } c_{Bv}^{\text{nom}} = \frac{\textcolor{gray}{\alpha_{B}} c_{B2}^{\text{nom}}}{\textcolor{gray}{\alpha_{A}} c_{A2}^{\text{nom}} + \textcolor{gray}{\alpha_{B}} c_{B2}^{\text{nom}} + \textcolor{gray}{\alpha_{C}} \left(1- c_{A2}^{\text{nom}} - c_{B2}^{\text{nom}} \right)}
\end{align}
As $\rho$ is always nonzero, we get $c_{B1}$ as:
\begin{align}
    c_{B1}(\rho) = \frac{c_{A1}(\rho) c_{B2}^{\text{nom}} - c_{A1}(\rho) c_{Bv}^{\text{nom}} + c_{A2}^{\text{nom}} c_{Bv}^{\text{nom}} - c_{Av}^{\text{nom}} c_{B2}^{\text{nom}}}{c_{A2}^{\text{nom}} - c_{Av}^{\text{nom}}}
\end{align}

From the third output component, we get the condition $\dot{\xi}_3 = \dot{T}_{2} = 0$ and thus solve Equation~(\ref{eq:DRMIMOSI_dT_2}) for $Q_2$:
\begin{align}
    \label{eq:DRMIMOSI_Q_2}
    Q_2(\rho,T_1) =-\textcolor{gray}{\varrho_{F}C_p} ( \rho +  F_{B}(\rho))(T_1- T_{2}^{\text{nom}})  + \textcolor{gray}{\Delta H_{V}} \rho .
\end{align}
The input $Q_2$ still depends on the temperature $T_1$ for which no expression has been derived up to this point.

From the fourth flat output component, we get $\dot{\xi}_4 = \dot{c}_{A1} = a_1^{\xi_4}\dot{\rho}$ (compare to Equation~(\ref{eq:DRMIMOSI_c_A2_rho})) and thus Equation~(\ref{eq:DRMIMOSI_dc_A1}) changes to:
    \begin{align}
        \label{eq:DRMIMOSI_for_T1}
        \textcolor{gray}{a_1^{\xi_4}}\dot{\rho} =& - c_{A1}(\rho) \textcolor{gray}{k_{1}} \textcolor{gray}{e}^{- \frac{\textcolor{gray}{E_{1}}}{\textcolor{gray}{R} T_{1}}} + \frac{\left(F_{B}(\rho) - F_{p}\right) \left(- c_{A1}(\rho) + c_{A2}^{\text{nom}}\right)}{\textcolor{gray}{V_{1}}} 
        \\ \nonumber &+  \frac{\left(F_{p} + \rho\right) \left(c_{A0} - c_{A1}(\rho)\right)}{\textcolor{gray}{V_{1}}}
    \end{align}
Equation~(\ref{eq:DRMIMOSI_for_T1}) is solved numerically for $T_1$. However, first, an expression for the input $F_p$ is needed. This expression is derived from the condition that according to our operating strategy the second derivative of the second output is zero, i.e., $\xi_2^{(2)} = c_{B2}^{(2)} = 0$. From this condition, we receive
    \begin{align}
        \label{eq:DRMIMOSI_F_p}
        F_p(\rho,T_1) = \Psi_{F_p}(\rho,T_1),
    \end{align}
with a nonlinear function $\Psi_{F_p}(\rho,T_1)$, which is not stated here to preserve readability. After inserting Equation~(\ref{eq:DRMIMOSI_F_p}) into Equation~(\ref{eq:DRMIMOSI_for_T1}), Equation~(\ref{eq:DRMIMOSI_for_T1}) implicitly defines $T_1$ as a function of $\rho$ and $\dot{\rho}$, $T_1(\rho, \dot{\rho})$. 
As discussed in the main paper, we cannot solve for $T_1$ analytically, but the solution received numerically is unique. 
With $T_1(\rho,\dot{\rho})$, we can calculate $Q_2$ as $Q_2(\rho, \dot{\rho})$ (compare to Equation~(\ref{eq:DRMIMOSI_Q_2})) and $F_p$ as $F_p(\rho, \dot{\rho})$ (compare to Equation~(\ref{eq:DRMIMOSI_F_p})).

The only variable that is missing is the reactor heating $Q_1$. 
To get an expression for $Q_1$, Equation~(\ref{eq:DRMIMOSI_for_T1}) is differentiated with respect to time.
Before doing so, we replace $F_p$ with $\Psi_{F_p}(\rho,T_1)$ (compare to Equation~(\ref{eq:DRMIMOSI_F_p})) and insert the functions for $F_B$ and $c_{A,1}$, which depend on $\rho$.
After bringing everything to the left hand side, Equation~(\ref{eq:DRMIMOSI_for_T1}) changes to
    \begin{align}
        0 = \Psi(\rho,\dot{\rho},T_1).
    \end{align}
Building the total differential with respect to time gives
    \begin{align}
        \label{eq:DRMIMOSI_Q_1}
        0 =& \frac{\partial \Psi(\rho,\dot{\rho},T_1)}{\partial \rho}\dot{\rho} + \frac{\partial \Psi(\rho,\dot{\rho},T_1)}{\partial \dot{\rho}}\nu 
        + \frac{\partial \Psi(\rho,\dot{\rho},T_1)}{\partial T_1}\Bigg[\frac{\rho + F_p}{\textcolor{gray}{V_1}}(\textcolor{gray}{T_{0}} - T_1) + 
        \\ \nonumber
        &\frac{F_B - F_p}{\textcolor{gray}{V_1}}(T_2 - T_1) - \frac{\textcolor{gray}{\Delta H_1}}{\textcolor{gray}{C_p}} \textcolor{gray}{k_1} c_{A1}\textcolor{gray}{e}^{-\frac{\textcolor{gray}{E_1}}{\textcolor{gray}{R}T_1}}- \frac{\textcolor{gray}{\Delta H_2}}{\textcolor{gray}{C_p}}\textcolor{gray}{k_2} c_{B1}\textcolor{gray}{e}^{-\frac{\textcolor{gray}{E_2}}{\textcolor{gray}{R}T_1}} + \frac{Q_1}{\textcolor{gray}{\varrho_F C_p V_1}} \Bigg ]
        \\ \nonumber
        =& \Psi_{0}(\rho,\dot{\rho},T_1) + \Psi_{\nu}(\rho,\dot{\rho},T_1)\nu + \Psi_{Q_1}(\rho,\dot{\rho},T_1)Q_1
        \\ \nonumber 
        &\text{with } \Psi_{\nu}(\rho,\dot{\rho},T_1) = \frac{\partial \Psi(\rho,\dot{\rho},T_1)}{\partial \dot{\rho}}\text{,  }\Psi_{Q_1}(\rho,\dot{\rho},T_1) =  \frac{\partial \Psi(\rho,\dot{\rho},T_1)}{\partial T_1}\frac{\textcolor{gray}{1}}{\textcolor{gray}{\varrho_F C_p V_1}},
        \\ \nonumber 
        &\text{and } \Psi_{0}(\rho,\dot{\rho},T_1) \text{ as the remaining part. }
    \end{align}
As $Q_1$ enters linearly, we can easily solve Equation~(\ref{eq:DRMIMOSI_Q_1}) for $Q_1$.
Moreover, the ramping degree of freedom $\nu$ enters linearly as well, which allows to solve Equation~(\ref{eq:DRMIMOSI_Q_1}) for $\nu$ and calculate the bounds of $\nu$ from the bounds of $Q_1$ (compare to Figure~ \ref{fig:DRMIMOSI_bounds_nu}).

At this point, all states and inputs can be calculated as functions of the production rate $\rho$ and its first two derivatives $\dot{\rho}$, and $\nu$.

Note that we do not derive expressions for the input derivatives $\dot{F}_B$, $F_B^{(2)}$, and $\dot{F}_p$ that are shown in Table~\ref{tab:DRMIMOSI_case_sudy_incidence}.
Instead, we replace $F_B$ with $F_B(\rho)$ (Equation~(\ref{eq:DRMIMOSI_F_B_of_rho})), and $F_p$ with $F_p(\rho,T_1)$  (Equation~(\ref{eq:DRMIMOSI_F_p})) in the equations where they appear and subsequently differentiate the complete equations with respect to time.

\subsection{Calculating the limits for $\dot{\rho}$}
As discussed in the previous section, the limits for $\nu$ can be calculated in a straightforward way from Equation~(\ref{eq:DRMIMOSI_Q_1}).
In this section, we calculate the limits for $\dot{\rho}$ resulting from the bounds of $T_1$, $Q_2$, $F_p$ (Figure 6 in the main paper).
First, the bounds are calculated as a function of $T_1$ by dividing Equation~(\ref{eq:DRMIMOSI_for_T1}) by $\textcolor{gray}{a_1^{\xi_4}}$, giving $\dot{\rho}$ as 
    \begin{align}
        \label{eq:DRMIMOSI_theata_T1}
        \dot{\rho} = \theta_{T_1}(\rho,T_1),
    \end{align}
with a nonlinear function $\theta_{T_1}$.
Note that while Equation~(\ref{eq:DRMIMOSI_for_T1}) cannot be solved for $T_1$ analytically, Equation~(\ref{eq:DRMIMOSI_for_T1}) can be solved analytically for $\dot{\rho}$.
Accordingly, inserting the bounds of $T_1$ gives the limits of $\dot{\rho}$ shown in Figure~6 in the main paper.

In the case of $Q_2$, which is calculated as $Q_2(\rho,T_1)$ from Equation~(\ref{eq:DRMIMOSI_Q_2}), Equation~(\ref{eq:DRMIMOSI_Q_2}) cannot be solved for $\dot{\rho}$ as no analytic expression was derived for $T_1$ that depends on $\rho$ and $\dot{\rho}$.
Instead, Equation~(\ref{eq:DRMIMOSI_Q_2}) is solved for $T_1$ as $T_1(Q_2,\rho)$, which is straightforward as $T_1$ enters Equation~(\ref{eq:DRMIMOSI_Q_2}) linearly.
Inserting $T_1(Q_2,\rho)$ into Equation~(\ref{eq:DRMIMOSI_theata_T1}) then gives
    \begin{align}
        \label{eq:DRMIMOSI_theata_Q_2}
        \dot{\rho} = \theta_{T_1}(\rho,T_1(Q_2,\rho)),
    \end{align}
which allows to calculate the limits of $\dot{\rho}$ from the bounds of $Q_2$ as shown in Figure 6 in the main paper.

For the purge stream $F_p$, we follow the same strategy to first calculate $T_1$ and then $\dot{\rho}$.
However, Equation~(\ref{eq:DRMIMOSI_F_p}) cannot be solved for $T_1$ as function of $F_p$ and $\rho$ analytically due to the exponential functions (see discussion around Equation~(11) in main paper).
Consequently, Equation~(\ref{eq:DRMIMOSI_F_p}) is solved numerically for $T_1$ as function of the bounds on $F_p$.
This value of $T_1$ is then used to calculate the limits of $\dot{\rho}$ shown in Figure 6 in the main paper.

\subsection{Ramping limits on $\nu$}
\begin{figure*}[h]
\centering
 \includegraphics[scale=0.8]{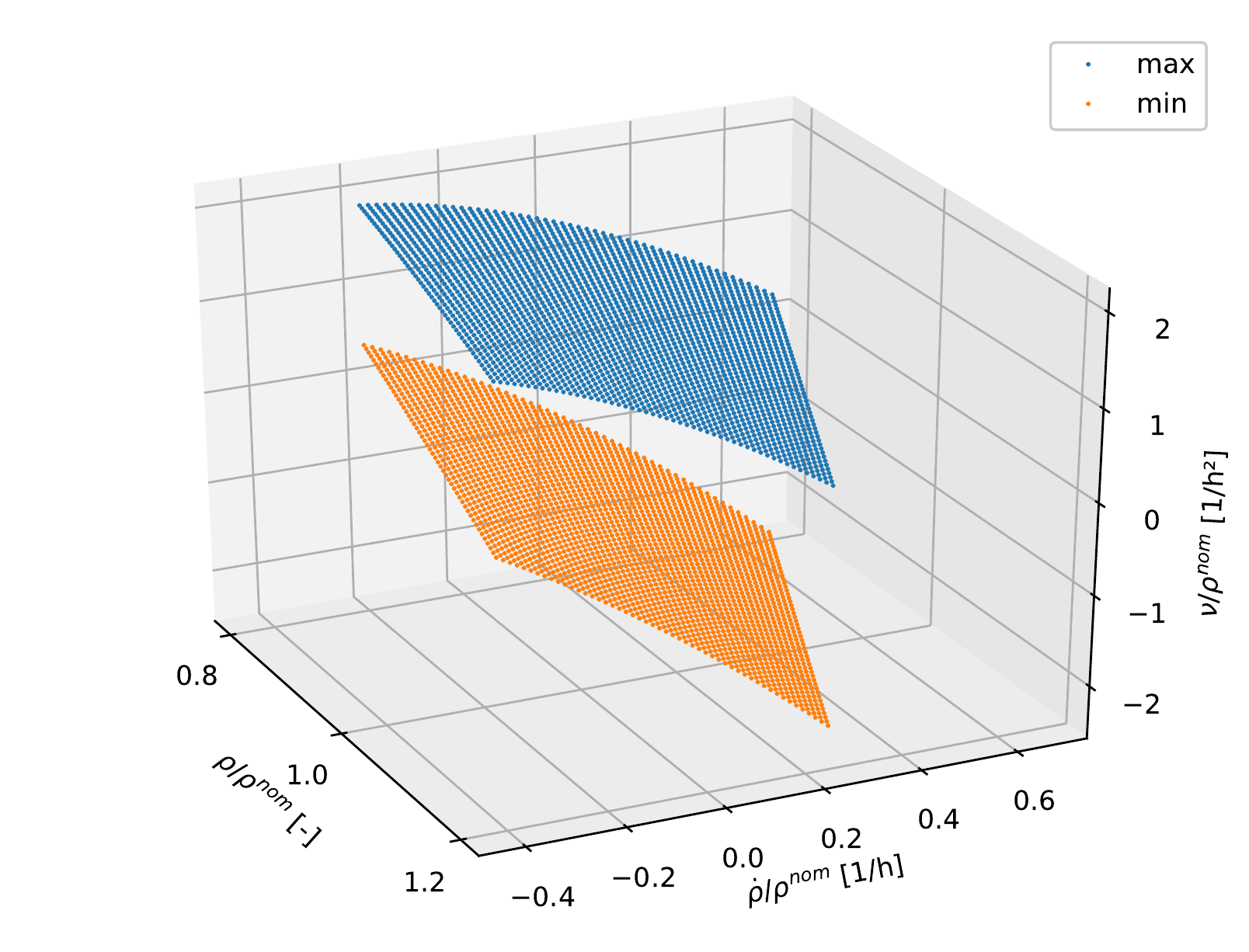}
  \caption{True limits on the ramping degree of freedom $\nu$ dependent on production rate $\rho$ and its first time derivative $\dot{\rho}$. All values are normalized with the nominal production rate $\rho^{\text{nom}}$}
   \label{fig:DRMIMOSI_bounds_nu}
\end{figure*}

In this section, the derivation of the ramping limits on the ramping degree of freedom $\nu$ is given in detail.
First, the space given by the limits $\rho^{\text{min}}$, $\rho^{\text{max}}$ and $\dot{\rho}^{\text{min}}(\rho)$, $\dot{\rho}^{\text{max}}(\rho)$ is sampled using 51$\times$51 points. 
For every point, the true limits $\nu_{\text{true}}^{\text{min}}(\rho, \dot{\rho})$, $\nu_{\text{true}}^{\text{max}}(\rho, \dot{\rho})$ are calculated from $Q_1^{\text{min}}$, $Q_1^{\text{max}}$ (compare to Equation~(\ref{eq:DRMIMOSI_Q_1})) (Figure \ref{fig:DRMIMOSI_bounds_nu}).
For a conservative approximation of these limits, we use piecewise affine functions. We distinguish between production rates $\rho$ above and below the nominal rate $\rho^{\text{nom}}$, and first derivatives $\dot{\rho}$ above and below zero.
Thus, the piecewise affine functions have 4 segments and need 2 binary variables in the optimization: $z_{\rho}$ indicating if $\rho$ is above or below the nominal rate $\rho^{\text{nom}}$, and $z_{\dot{\rho}}$ indicating if $\dot{\rho}$ is above or below zero.
In every segment $s$, the limits $\nu_{\text{PWA}}^{\text{min}}(\rho, \dot{\rho})$, $\nu_{\text{PWA}}^{\text{max}}(\rho, \dot{\rho})$ are given as
\begin{align}
    &\nu_{\text{PWA},s}^{\text{min}}(\rho, \dot{\rho}) = a_{0,s}^{\text{min}} + a_{\rho,s}^{\text{min}} \rho + a_{\dot{\rho},s}^{\text{min}}\dot{\rho} , \\
    &\nu_{\text{PWA},s}^{\text{max}}(\rho, \dot{\rho}) = a_{0,s}^{\text{max}} + a_{\rho,s}^{\text{max}} \rho + a_{\dot{\rho},s}^{\text{max}}\dot{\rho} ,
\end{align}
with parameters $a_{0,s}^{\text{min}}$, $a_{\rho,s}^{\text{min}}$, $a_{\dot{\rho},s}^{\text{min}}$, $a_{0,s}^{\text{max}}$, $a_{\rho,s}^{\text{max}}$, $a_{\dot{\rho},s}^{\text{max}}$.
The parameters are calculated in an optimization, which minimizes the difference to the true bounds while ensuring conservativeness such that $\nu_{\text{true}}^{\text{min}}(\rho, \dot{\rho}) \leq \nu_{\text{PWA},s}^{\text{min}}(\rho, \dot{\rho}) < \nu_{\text{PWA},s}^{\text{max}}(\rho, \dot{\rho}) \leq \nu_{\text{true}}^{\text{max}}(\rho, \dot{\rho})$ holds for every point $(\rho, \dot{\rho})$.

To evaluate the quality of the bounds, we calculate the coverage $C(\rho, \dot{\rho})$:
\begin{align}
    C(\rho, \dot{\rho}) = \frac{\nu_{\text{PWA}}^{\text{max}}(\rho, \dot{\rho}) - \nu_{\text{PWA}}^{\text{min}}(\rho, \dot{\rho})}{\nu_{\text{true}}^{\text{max}}(\rho, \dot{\rho}) - \nu_{\text{true}}^{\text{min}}(\rho, \dot{\rho})},
\end{align}
which is the distance between the chosen bounds in relation to the distance between the true bounds.
The arithmetic mean of the coverage is 95~\%.
Accordingly, the chosen piecewise affine limits capture approximately 95~\% of the feasible region.
Therefore, we do not refine the approximation by introducing more segments.

\subsection{Multi-energy system and optimization problem for demand response application}
In this section, we describe the multi-energy system and the optimization problem in the second study, which considers a demand response application.
Our aim for this case study is to construct an illustrative example that demonstrates the application and capabilities of our new method.
The considered multi-energy system is based on the third day of the benchmark problem given by \cite{SusanneSass.2020}. 
Note that the data are publicly available at \cite{HECI}. 
The given heat and electricity demands are multiplied by 3 as we increase the number of units from 2 to 6.

Instead of one combined heat and power plant (CHP) and one boiler (B) as in \cite{SusanneSass.2020}, we use 4 CHPs and 2 boilers.
All 4 CHPs have a nominal thermal output power of 450 kW, and the 2 boilers have 530 kW.
The nominal efficiencies are given in Table \ref{tab:DRMIMOSI_effs}.

\begin{table}[h!]
    \centering
    \caption{Nominal efficiencies of energy-system components.}
    \label{tab:DRMIMOSI_effs}
    \begin{tabular}{lrr} \\
    unit & thermal efficiency & electric efficiency  \\
     \hline
    B1 & 79.2\%  & - \\
    B1 & 80.8\%  & - \\
    CHP1 & 48.3\% & 37.7\%\\
    CHP2 & 49.3\%  & 38.4\% \\
    CHP3 & 50.3\%  & 39.2\% \\
    CHP4 & 51.3\%  & 39.9\% \\
    \end{tabular}
\end{table}

The minimum part-load fraction is 50\% for the CHPs and 20\% for the boilers \citep{SusanneSass.2020}.
Following \cite{Voll.2014}, the part-load efficiency curves are discretized with one affine element to achieve an accurate discretization.

The electricity price used is the German day-ahead price from April 2nd, 2021 \citep{SMARD}.
This price curve is chosen because it features times with low prices that make the operation of the boiler favorable and times of high prices which make the operation of the CHPs favorable. 

For the optimization problem, which schedules the flexible operation of the reactor-separator process, a process energy demand model is needed (Pc).
As the dynamic ramping constraint is second order, the heat demand of the process is modeled as a function of $\rho$, $\dot{\rho}$, and $\nu$, i.e.,
\begin{align}
    Q_{\text{dem,heat}}^{\text{process}}(\rho, \dot{\rho}, \nu).
\end{align}
To fit the heat demand model, we sample the operating region given by the limits $\rho^{\text{min}}$, $\rho^{\text{max}}$, $\dot{\rho}^{\text{min}}(\rho)$, $\dot{\rho}^{\text{max}}(\rho)$, and $\nu^{\text{min}}(\rho, \dot{\rho})$, $\nu^{\text{max}}(\rho, \dot{\rho})$ using 11$\times$11$\times$11 points.
First, a purely linear function is used for the heat demand and the parameters are fitted using linear regression.
Because the average absolute approximation error is 10\% of the nominal heat demand, we try a piecewise-affine (PWA) model to increase the accuracy.
For this new PWA model, the operating region is split into segments in which affine functions are used.
Here, the operating region is split into four segments by making two distinctions: If $\dot{\rho}$ is greater or smaller than zero and if $\nu$ is greater or smaller than zero.
To implement the PWA model in the optimization, two binary variables are needed: $z_{\dot{\rho}}$ indicates whether $\dot{\rho}$ is positive and $z_{\nu}$ indicates whether $\nu$ is positive.
No distinction is made for $\rho$ as we find that the heat demand is roughly linear with $\rho$.
The PWA model reduces the mean approximation error to 5\% of the nominal heat demand.

The optimization problem is discretized using orthogonal collocation with 2 elements per hour and 3 collocation points per element.
The binary variables are discretized with one-hour timestep to match the hourly changing prices and demands.
There are 9 binary variables per timestep: 4 on/off binaries for the CHPs, 2 on/off binaries for the boilers, 2 binaries $z_{\rho}$, $z_{\dot{\rho}}$ for the piecewise affine approximation of the ramping limits $\nu_{\text{PWA}}^{\text{min}}(\rho, \dot{\rho})$, $\nu_{\text{PWA}}^{\text{max}}(\rho, \dot{\rho})$, and 1 additional binary $z_{\nu}$ for the PWA heat demand model. 
Consequently, the optimization problem has 9$\times$24=216 binary variables for the 24~h time horizon.
The ramping degree of freedom $\nu$ is discretized to be piecewise linear for every hour.

\bibliographystyle{apalike}
\renewcommand{\refname}{Bibliography} 
\bibliography{literature.bib}